\documentclass[journal,twoside,web]{ieeecolor}
\usepackage{generic}
\usepackage{cite}
\usepackage{amsmath,amssymb,amsfonts}
\usepackage{algorithmic}
\usepackage{graphicx}
\usepackage{textcomp}
\def\BibTeX{{\rm B\kern-.05em{\sc i\kern-.025em b}\kern-.08em
    T\kern-.1667em\lower.7ex\hbox{E}\kern-.125emX}}
\usepackage{amscd}
\usepackage{dsfont}
\usepackage{fancyhdr}
\usepackage{mathtools}

\usepackage{cite}
\usepackage[ruled]{algorithm2e}
\usepackage{ifthen}
\usepackage{algorithmic}
\newtheorem{theorem}{Theorem}
\newtheorem{lemma}{Lemma}

\newtheorem{proposition}{Proposition}
\newtheorem{corollary}{Corollary}

\usepackage{empheq}

\def\eq#1{\begin{equation}#1\end{equation}}
\newcommand{\R}{{\rm I\!R}}
\newcommand{\bbb}{\mathbb}

\newcommand{\1}{\mathbf{1}}

\def\##1\#{\begin{align}#1\end{align}}
\def\$#1\${\begin{align*}#1\end{align*}}

\def\qed{ \rule{.08in}{.08in}}

\newcommand{\dfb}{\stackrel{\Delta}{=}}

\usepackage{color}
\usepackage{hyperref}

\begin{document}
\title{A Class of Optimal Directed Graphs for Network Synchronization}
\author{Susie Lu, \IEEEmembership{Student Member, IEEE}, John Urschel, Ji Liu, \IEEEmembership{Member, IEEE}
\thanks{S.~Lu and J.~Urschel are with the Department of Mathematics at the Massachusetts Institute of Technology.
J.~Liu is with the Department of Electrical and Computer Engineering at Stony Brook University.
Email addresses: \href{mailto:susielu@mit.edu}{susielu@mit.edu}, \href{mailto:urschel@mit.edu}{urschel@mit.edu}, \href{mailto:ji.liu@stonybrook.edu}{ji.liu@stonybrook.edu}.}
}

\maketitle

\begin{abstract}
In a paper by Nishikawa and Motter, a quantity called the normalized spread of the Laplacian eigenvalues is used to measure the synchronizability of certain network dynamics. Through simulations, and without theoretical validation, it is conjectured that among all simple directed graphs with a fixed number of vertices and arcs, the optimal value of this quantity is achieved if the Laplacian spectrum satisfies a specific pattern. This paper proves this conjecture and further shows that the conjectured spectral condition is not only sufficient but also necessary. Moreover, the paper proves that the optimal Laplacian spectrum is always achievable by a class of almost regular directed graphs, which can be constructed through an inductive algorithm. 
\end{abstract}

\section{The Conjecture}

Over the past two decades, synchronization in complex networks has attracted considerable attention for its crucial role in fields including biology, climatology, ecology, sociology, and technology \cite{physicsreport}. A typical class of network synchronization dynamics can be described as
\eq{\textstyle \dot{x}_i = F(x_i)+\frac{\varepsilon}{d} \sum_{j=1}^n a_{ij} \big(H(x_j) - H(x_i)\big),\label{eq:dynamics}}
where $x_i$ is the state vector of the $i$th dynamical agent in a network of $n$ agents, $F$ represents the dynamics of each agent when isolated, $H(x_j)$ denotes the signal that the $j$th agent sends to its neighboring agents, $a_{ij}$ is the $ij$th entry of the adjacency matrix of the underlying simple directed graph, and $\varepsilon/d$ represents the coupling strength $\varepsilon$ normalized by the average coupling strength per vertex $d = \frac{1}{n} \sum_{i=1}^n \sum_{j=1}^n a_{ij}$.
The network synchronization problem is to derive conditions under which all $n$ agents' states converge to the same stable state. 
More details on the above network synchronization dynamics can be found in \cite{Nishikawa10}. 

In the case when the underlying graph is undirected, the synchronizability of the network dynamics is measured by the eigenratio, which is defined as the ratio of the largest eigenvalue to the smallest nonzero eigenvalue of the Laplacian matrix \cite[Section 4.1.2]{physicsreport}. For directed graphs, however, there is no standard index for measuring synchronizability \cite[Page 115]{physicsreport}.
It is worth emphasizing that the network synchronization dynamics described in \cite[Equation (54)]{physicsreport} is mathematically equivalent to \eqref{eq:dynamics}, although they use slightly different notation. 

A directed graph is called simple if it does not have any self-arcs and multiple directed edges with the same tail and head vertices.
For any simple
directed graph $\bbb G$ with $n$ vertices, we use $D$ and $A$ to denote its in-degree matrix and adjacency matrix, respectively. Specifically, $D$ is an $n\times n$ diagonal matrix whose $i$th diagonal entry equals the in-degree of vertex $i$, and $A$ is an $n\times n$ matrix whose $ij$th entry equals $1$ if $(j,i)$ is an arc (or a directed edge) in $\bbb G$ and otherwise equals $0$. The Laplacian matrix of $\bbb G$ is defined as $L=D-A$. It is easy to see that a Laplacian matrix $L$ always has an eigenvalue at 0 since all its row sums equal 0. 
In the special case when $\bbb G$ is a simple undirected graph, each undirected edge between two vertices $i$ and $j$ can be equivalently replaced by a pair of directed edges $(i,j)$ and $(j,i)$; then $L$ is a symmetric matrix and thus has a real spectrum. It is well known that in this case $L$ is positive-semidefinite, its smallest eigenvalue equals 0, and its second smallest eigenvalue is positive if and only if $\bbb G$ is connected \cite{Fiedler73}. 
For directed graphs, $L$ may have complex eigenvalues. 
Let $\lambda_1,\lambda_2,\ldots,\lambda_n$ denote all $n$ eigenvalues of $L$, with $\lambda_1=0$ and $\lambda_2,\ldots,\lambda_n$ possibly complex. 
Define
\eq{\textstyle\sigma^2 \dfb \frac{1}{n-1} \sum_{i=2}^n |\lambda_i - \bar\lambda|^2, \;\; {\rm where} \;\; \bar\lambda = \frac{1}{n-1} \sum_{i=2}^n \lambda_i,\label{eq:sigma}}
which is a normalized deviation of possibly nonzero eigenvalues. This quantity is called\footnote{
The definition of the normalized spread of the eigenvalues in \cite[Page 10343]{Nishikawa10} includes an additional $d^2$ term in the denominator; in this paper, we adopt a simplified definition without the $d^2$ term, as its omission does not affect the conjecture in \cite{Nishikawa10} for unweighted graphs.
} the normalized spread of the eigenvalues in \cite{Nishikawa10} to measure the synchronizability 
of \eqref{eq:dynamics}.
It is claimed and validated by simulations that the smaller the value of $\sigma^2n^2/m^2$, the more synchronizable the network will generally be, where $m$ denotes the number of directed edges in $\bbb G$.  Note that since the sum of all $n$ agents' in-degrees equals $m$ and the sum of all eigenvalues of a matrix equals the trace of the matrix, it follows that $\bar\lambda = m/(n-1)$ which is always a real number.
It is clear that, for fixed values of $n$ and $m$, the smallest possible $\sigma$ corresponds to the optimal graph(s) for network synchronization.
The following conjecture was proposed by Nishikawa and Motter in \cite[Page 10343]{Nishikawa10}.

\vspace{.05in}

{\bf Conjecture:}
Among all simple directed graphs with $n$ vertices and $m$ arcs, the minimum possible value of $\sigma^2$ is
\eq{\textstyle \sigma^2_{\min} = \frac{1}{(n-1)^2} \big[m-(n-1)\kappa\big]\big[(n-1)(\kappa+1)-m\big], 
\label{eq:spread}}
which is achieved if the Laplacian spectrum is 
\eq{
0, \underbrace{\kappa,\; \ldots, \;\kappa,}_{(n-1)(\kappa+1)-m} \; \underbrace{\kappa+1, \ldots, \kappa+1}_{m-(n-1)\kappa}, \label{eq:spectrum}
}
where $\kappa \dfb \lfloor \frac{m}{n-1}\rfloor$. Here $\lfloor \cdot \rfloor$ denotes the floor function.



Note that $\kappa$ is the unique integer such that $(n-1)\kappa\le m < (n-1)(\kappa+1)$. Defining $q_\kappa=(n-1)\kappa$ and $q_{\kappa+1}=(n-1)(\kappa+1)$, it follows that\footnote{There is a typo in \cite{Nishikawa10} which states $q_\kappa \le m \le q_{\kappa+1}$.} $q_\kappa \le m < q_{\kappa+1}$ and $\sigma_{\min}^2 = (m-q_\kappa)(q_{\kappa+1}-m)/(n-1)^2$, which is consistent with the expressions in \cite{Nishikawa10}. 
It was implicitly assumed in \cite{Nishikawa10} that $m\ge n-1$ for the conjecture, as indicated by its Fig. 2A where $m$ ranges from $n-1$ to $n(n-1)$. This assumption is natural, as network synchronization requires connectivity, and $n-1$ is the minimal number of arcs needed to guarantee a connected network. We will show that the conjecture holds even when $m<n-1$.  

It is worth emphasizing that the conjecture itself is independent of network synchronization dynamics, even though it was proposed as an optimal synchronization condition.

There is another popular way to define a Laplacian matrix of a directed graph based on out-degree \cite{wiki}. Specifically, letting $D_{{\rm out}}$ denote the out-degree matrix of $\bbb G$, the corresponding Laplacian matrix is denoted and defined as $L_{{\rm out}}=D_{{\rm out}}-A'$. It is straightforward to verify that the in-degree Laplacian matrix of a directed graph $\bbb G$ equals the out-degree Laplacian matrix of its transpose graph
$\bbb G'$, where the transpose of a directed graph is the directed graph with the same vertex set obtained by reversing all its directed edges. In other words, $L(\bbb G)= L_{{\rm out}}(\bbb G')$. Since the set of all possible simple directed graphs with $n$ vertices and $m$ arcs is invariant under the graph transpose operation, the conjecture remains unaffected regardless of whether the Laplacian matrix is defined based on in-degree or out-degree. The only resulting difference is that any optimal graph with $\sigma_{\min}$ needs to be correspondingly transposed if the definition is changed from in-degree based to out-degree based. 

The conjecture was validated in \cite{Nishikawa10} only for small-sized graphs with $n \le 6$ through simulations; no theoretical validation was provided therein. Indeed, to the best of our knowledge, the conjecture has never been studied from a theoretical perspective. It even remains unclear whether the conjectured Laplacian spectrum \eqref{eq:spectrum} exists for a given fixed number of vertices $n$ and arcs $m$.
This is exactly what this paper aims to address.

This paper proves the conjecture by establishing the following theorem, which states a stronger version of the conjecture.

\vspace{.05in}

\begin{theorem}\label{thm:main}
Among all simple directed graphs with $n$ vertices and $m$ arcs, the minimum possible value of $\sigma^2$ is \eqref{eq:spread}, which is achieved if, and only if, the Laplacian spectrum is~\eqref{eq:spectrum}.
\end{theorem}

\vspace{.05in}

The theorem establishes that the conjectured minimal value of $\sigma^2$ is attainable and that the conjectured spectral condition is not only sufficient but also necessary for optimal synchronizability.
The paper further shows that, for any feasible pair of $n$ and $m$, the conjectured Laplacian spectrum can always be achieved by a class of ``almost regular'' directed graphs. Consequently, all these graphs are optimal graphs for network synchronization. Moreover, for any fixed number of vertices $n$, these graphs can be generated through an inductive construction algorithm for each possible number of arcs $m$.
The algorithm was presented in a preliminary conference version \cite{cdc25}, in which the conjecture was proved only for a few special cases. This paper provides a complete analysis of the algorithm and its generated graphs, and proves the conjecture in full generality, which were not included in \cite{cdc25}.

\section{Proof of the Conjecture}

To prove the conjecture, and in fact more strongly, to prove Theorem~\ref{thm:main}, the analysis consists of two main logical steps. The first is to establish that the conjectured Laplacian spectrum \eqref{eq:spectrum} is actually a necessary and sufficient condition for achieving $\sigma^2_{\min}$ in \eqref{eq:spread}, the minimal normalized spread of the Laplacian eigenvalues (cf. Corollary \ref{coro:pnas}); the second is to show that this conjectured optimal condition is always realizable (cf. Theorem \ref{th:pnas-spectrum}).
These two steps will be addressed in the following two subsections, respectively.

\subsection{Necessary and Sufficient Spectral Condition}

The necessity and sufficiency of the Laplacian spectral condition \eqref{eq:spectrum} are a consequence of the following theorem.

\vspace{.05in}

\begin{theorem}\label{th:generalized}
For any integer-coefficient monic polynomial $p(x)$ of positive degree $k$ with complex roots $r_1, \ldots, r_k$ whose sum $\sum_{i=1}^k r_i$ equals $\ell$, 
\eq{
\textstyle \frac{1}{k} \sum_{i=1}^k \big| r_i - \frac{\ell}{k} \big|^2 \ge \big(\frac{\ell}{k} - \lfloor \frac{\ell}{k} \rfloor \big) \big( \lceil \frac{\ell}{k} \rceil - \frac{\ell}{k} \big),
\label{eq:generalized}}
and equality holds if, and only if, each root equals $\lfloor \frac{\ell}{k} \rfloor$ or $\lfloor \frac{\ell}{k} \rfloor + 1$, where $\lceil \cdot \rceil$ denotes the ceiling function.
\end{theorem}

\vspace{.05in}

By Vi\`ete's formulas, the sum of the roots is the negative of the coefficient of $x^{k-1}$ in $p(x)$, and hence $\ell$ is an integer.
Also, it is worth emphasizing that when all roots are integers and their difference is at most one, for example when each root equals $\lfloor \ell/k \rfloor$ or $\lfloor \ell/k \rfloor + 1$, and their sum is a given integer, the multiset of roots is uniquely determined.

For any simple directed graph with $n$ vertices and $m$ arcs, let $p_L(x)=\det (xI-L)$ be the characteristic polynomial of its Laplacian matrix $L$, whose eigenvalues are $\lambda_1,\lambda_2,\ldots,\lambda_n$ with $\lambda_1=0$. Then, $q(x)=p_L(x)/x$ is an integer-coefficient monic polynomial of degree $n-1$ whose complex roots are $\lambda_2,\ldots,\lambda_n$. 
Since their sum $\sum_{i=2}^n \lambda_i$ equals the number of arcs $m$, which is an integer, Theorem~\ref{th:generalized} applies here with $k=n-1>0$ and $\ell = m$. 
First, note that $\sigma^2 = \frac{1}{n-1} \sum_{i=2}^n |\lambda_i - \bar\lambda|^2 = \frac{1}{n-1} \sum_{i=2}^n |\lambda_i - \frac{m}{n-1}|^2$, which equals the left hand side of \eqref{eq:generalized}.
Second, with $k = n - 1$ and $\ell = m$, the right hand side of \eqref{eq:generalized} can be written as 
$$
\textstyle \frac{1}{(n-1)^2} \big(m-(n-1)\lfloor \frac{m}{n-1} \rfloor\big)\big((n-1)\lceil \frac{m}{n-1} \rceil-m\big). 
$$
Note that if $m$ is a multiple of $n-1$, then $m-(n-1)\lfloor \frac{m}{n-1} \rfloor = 0$, and if $m$ is not a multiple of $n-1$, then $\lceil \frac{m}{n-1} \rceil = \lfloor \frac{m}{n-1} \rfloor +1$. In either case, the above expression coincides with that in \eqref{eq:spread}. 
Last, under the condition $\sum_{i=2}^n \lambda_i=m$, if each of $\lambda_2,\ldots,\lambda_n$ equals either $\lfloor \frac{m}{n-1} \rfloor$ or $\lfloor \frac{m}{n-1} \rfloor + 1$, then the Laplacian spectrum must be 
given by \eqref{eq:spectrum}, and clearly vice versa. 
We therefore have proved the following corollary.

\vspace{.05in}

\begin{corollary}\label{coro:pnas}
For any simple directed graph with $n$ vertices and $m$ arcs, 
$\sigma^2 \ge \frac{1}{(n-1)^2} [m-(n-1)\kappa][(n-1)(\kappa+1)-m]$, 
and equality holds 
if, and only if, the Laplacian spectrum is~\eqref{eq:spectrum}.
\end{corollary}

\vspace{.05in}


{\bf Proof of Theorem \ref{th:generalized}:}
Let $a = \lfloor \frac{\ell}{k} \rfloor$ and $b = \frac{\ell}{k} - \lfloor \frac{\ell}{k} \rfloor$. Express $p(x)$ as $(x-a)^\alpha (x-a-1)^\beta q(x)$, where $\alpha$ and $\beta$ are nonnegative integers, and $q(x)$ is a monic integer-coefficient polynomial with no roots at $a$ or $a+1$.
It is clear that $q(x)$ is of order $\gamma = k - \alpha - \beta \ge 0$, and we write its roots as $s_1, \ldots, s_{\gamma}$. Then, 
\begin{align}
    & \textstyle\sum_{i=1}^k | r_i - \frac{\ell}{k}|^2 = \textstyle\sum_{i=1}^k | r_i - (a+b)|^2 \nonumber\\
    &= \alpha b^2 + \beta (1-b)^2 + \textstyle\sum_{i=1}^\gamma |s_i - (a+b)|^2. \label{eq:sum}
\end{align}
Note that the right hand side of \eqref{eq:generalized}, $(\frac{\ell}{k} - \lfloor \frac{\ell}{k} \rfloor ) ( \lceil \frac{\ell}{k} \rceil - \frac{\ell}{k} )$, equals $b(1-b)$. 
To prove the theorem, it is equivalent to show that the expression in \eqref{eq:sum} is no smaller than $k b(1-b)$, and that equality holds if 
and only if each $r_i$ equals $a$ or $a+1$.

When $b = 0$, equivalently when $\ell$ is a multiple of $k$, the right hand side of inequality \eqref{eq:generalized} equals $0$, and in this case the theorem is clearly true. It is therefore assumed that $b > 0$ for the remainder of the proof.

We first consider the case when $\gamma = 0$, which implies that $\alpha + \beta = k$ and all roots of $p(x)$ are $a$ or $a+1$.
In this case, $a+b = \frac{\ell}{k} = \frac{1}{k} \sum_{i=1}^k r_i = \frac{1}{k}(\alpha a + \beta (a+1)) = a + \frac{\beta}{k}$, so $\beta =kb$, and thus $\alpha = k(1-b)$. 
From \eqref{eq:sum}, 
\begin{align*}
    \textstyle \sum_{i=1}^k | r_i - \frac{\ell}{k}|^2 = \alpha b^2 + \beta (1-b)^2
    = kb(1-b),
\end{align*}
and therefore equality holds in \eqref{eq:generalized} in this case.

We next consider the case when $\gamma$ is a positive integer.  
Note that $\sum_{i=1}^\gamma s_i = \sum_{i=1}^k r_i - \alpha a - \beta (a+1) = \ell - \alpha a - \beta (a+1)$ is an integer. 
We divide the analysis into two cases depending on whether 
$\frac{1}{\gamma}\sum_{i=1}^\gamma s_i$ lies in the interval $[a, a+1]$ or not. 

{\bf Case 1:} Suppose that $\frac{1}{\gamma} \sum_{i=1}^\gamma s_i \in [a, a+1]$.
Express the monic polynomial $q(x)$ as $\prod_{i=1}^\gamma (x - s_i)$.
From the inequality of arithmetic and geometric means, for any integer $z$ that is not a root of $q(x)$,
\eq{ \textstyle\sum_{i=1}^\gamma |s_i - z|^2 \ge \gamma \big(\prod_{i=1}^\gamma |s_i - z| \big)^{\frac{2}{\gamma}} = \gamma |q(z)|^{\frac{2}{\gamma}} \ge \gamma, \label{eq:cake}}
where we used the fact that $q(z)$ is a nonzero integer in the last inequality.
It is straightforward to verify that 
\begin{align*}
|s_i - (a+b)|^2 
=\;& b|s_i - (a+1)|^2 + (1-b) |s_i - a|^2 \\
&\!- b(1-b),
\end{align*}
which expresses $|s_i - (a+b)|^2$ in terms of $|s_i - a|^2$ and $|s_i - (a+1)|^2$.
Substituting this expression into \eqref{eq:sum} and applying inequality \eqref{eq:cake} yields
\begin{align}
\textstyle\sum_{i=1}^k \left| r_i - \frac{\ell}{k}\right|^2 =\;& \textstyle \alpha b^2 + \beta (1-b)^2 +  b \sum_{i=1}^\gamma |s_i - (a+1)|^2 \nonumber\\
&\! + \textstyle (1-b) \sum_{i=1}^\gamma |s_i - a|^2 - \gamma b(1-b) \nonumber\\
\ge\;&  \alpha b^2 + \beta (1-b)^2 + \gamma (1-b+b^2). \label{eq:case1bound}
\end{align}
We claim that, among all feasible triples $(\alpha, \beta, \gamma)$, the lower bound in \eqref{eq:case1bound} attains its smallest value only when $\alpha+\beta = k$, or equivalently $\gamma = 0$.
To prove the claim, we use $\tilde{\alpha}$ and $\tilde{\beta}$ to respectively denote the values of $\alpha$ and $\beta$ in the case $\gamma = 0$. From the preceding discussion, these values are unique, with $\tilde{\alpha} = k(1-b)$ and $\tilde{\beta} = kb$.
Recall that any triple $(\alpha, \beta, \gamma)$ satisfies the two equalities $\alpha+\beta+\gamma=k$ and $\alpha a + \beta (a+1) + \sum_{i=1}^\gamma s_i = \ell$.
In the case $\gamma = 0$, these two equalities simplify to $\tilde\alpha+\tilde\beta=k$ and $\tilde\alpha a + \tilde\beta (a+1) = \ell$.
These four equalities lead to the following linear equations in $(\tilde\alpha-\alpha)$ and $(\tilde\beta-\beta)$: 
\begin{align*}
    (\tilde\alpha-\alpha) a + (\tilde\beta-\beta) (a+1) &= \textstyle\sum_{i=1}^\gamma s_i, \\
    (\tilde\alpha-\alpha)  + (\tilde\beta-\beta)  &= \gamma.
\end{align*}
Solving the above linear equations yields
\begin{align*}
    \tilde\alpha-\alpha &= \gamma (a+1) - \textstyle\sum_{i=1}^\gamma s_i, \\
    \tilde\beta-\beta  &= \textstyle\sum_{i=1}^\gamma s_i - \gamma a.
\end{align*}
Since $\frac{1}{\gamma} \sum_{i=1}^\gamma s_i \in [a, a+1]$, it is easy to see that both $(\tilde\alpha-\alpha)$ and $(\tilde\beta-\beta)$ are nonnegative. Moreover, their sum equals $\gamma>0$, and thus at least one of them is positive. 
In other words, $\tilde{\alpha} \ge \alpha$ and $\tilde{\beta} \ge \beta$, with at least one of these inequalities being strict.
Substituting $\gamma = k - \alpha - \beta$, the lower bound in \eqref{eq:case1bound} can be written as the following function of $\alpha$ and $\beta$: 
$$\alpha (b-1) - \beta b + k(1-b+b^2).$$
Since $0 < b < 1$, the above function is strictly decreasing in $\alpha$ and $\beta$. 
It follows that the function achieves its minimum uniquely at $\alpha = \tilde{\alpha}$ and $\beta = \tilde{\beta}$, that is, when $\gamma = 0$. This proves the claim. 
Substituting $\tilde{\alpha} = k(1-b)$ and $\tilde{\beta} = kb$ into the function, its minimum equals $kb(1-b)$. 
From \eqref{eq:case1bound}, we conclude that in Case 1, $\sum_{i=1}^k | r_i - \frac{\ell}{k} |^2 > kb(1-b)$.

{\bf Case 2:} Suppose that $\frac{1}{\gamma} \sum_{i=1}^\gamma s_i \notin [a, a+1]$. Let $s = \frac{1}{\gamma} \sum_{i=1}^\gamma s_i$. Then, either $s<a$ or $s>a+1$. Note that
\begin{align*}
    &\textstyle\sum_{i=1}^\gamma |s_i - (a+b)|^2 \\
    =\; & 
    \textstyle\sum_{i=1}^\gamma |s_i|^2 - (a+b) \sum_{i=1}^\gamma (s_i + \bar s_i) + \gamma (a+b)^2 \\
    =\; & 
    \textstyle\sum_{i=1}^\gamma |s_i|^2 - 2s\gamma(a+b) + \gamma (a+b)^2 \\
    \overset{(a)}{\ge}\; &
    \textstyle \frac{1}{\gamma} |\sum_{i=1}^\gamma s_i |^2 - 2s\gamma(a+b) + \gamma (a+b)^2 \\
    =\; & s^2\gamma - 2s\gamma(a+b) + \gamma (a+b)^2
    = \gamma (s-a-b)^2,
\end{align*}
where in (a) we used the Cauchy--Schwarz inequality.
Plugging the above inequality into \eqref{eq:sum} yields
\begin{align}
    \textstyle\sum_{i=1}^k | r_i - \frac{\ell}{k}|^2 
    \ge  \alpha b^2 + \beta (1-b)^2 + \gamma (s-a-b)^2. \label{eq:case2bound}
\end{align}
We claim that $\alpha b^2 + \beta (1-b)^2 + \gamma (s-a-b)^2 > kb(1-b)$. To see this, let $c=s-a$, and then
$$\alpha b^2 + \beta (1-b)^2 + \gamma (s-a-b)^2 = kb^2 + \beta(1-2b) + \gamma c (c - 2b),$$
where we used the fact that $\alpha+\beta+\gamma = k$.
It follows that proving the claim is equivalent to proving
\eq{\beta(1-2b) + \gamma c (c - 2b) > kb(1-2b). \label{eq:pen}}
Recall the facts that $a+b=\frac{\ell}{k}$ and $\ell = \alpha a +\beta (a+1) + \sum_{i=1}^\gamma s_i = \alpha a +\beta (a+1) + \gamma s$. It follows that 
$kb = \ell - ka = \alpha a +\beta (a+1) + \gamma s - ka = \beta + \gamma s - \gamma a = \beta + \gamma c$.  
Thus, the above inequality \eqref{eq:pen} is equivalent to 
$\beta(1-2b) + \gamma c (c - 2b) > (\beta +\gamma c)(1-2b)$, which further simplifies in a straightforward manner to $c(c-1)>0$. 
We prove this inequality by considering the two possible cases for $s$. First, if $s<a$, then $c<0$, and thus $c(c-1)>0$. Second, if $s>a+1$, then $c>1$, and $c(c-1)>0$ also holds. 
This proves the claim. From \eqref{eq:case2bound}, $\sum_{i=1}^k | r_i - \frac{\ell}{k}|^2 > kb(1-b)$ in Case 2.

From the preceding discussion, we conclude that when $b > 0$, $\sum_{i=1}^k | r_i - \frac{\ell}{k}|^2$ is no smaller than $k b(1-b)$, and that equality holds if and only if $\gamma = 0$, that is, if and only if each $r_i$ is either $a$ or $a+1$.
This completes the proof of the theorem.
\hfill$\qed$

\subsection{Optimal Graphs}

In this subsection, we present an algorithm that, for any fixed number of vertices $n > 1$, inductively constructs a class of simple directed graphs, each having the Laplacian spectrum specified in \eqref{eq:spectrum} for every possible number of arcs $m\ge n-1$. 
Note that if $m < n-1$, the graph is disconnected, in which case network synchronization cannot be achieved.
To this end, we need the following concepts.

A vertex $i$ in a directed graph $\bbb G$ is called a root of $\bbb G$ if for each other vertex $j$ of $\bbb G$, there is a directed path from $i$ to $j$. We say that $\bbb G$ is rooted at vertex $i$ if $i$ is in fact a root, and that $\bbb G$ is rooted if it possesses at least one root. In other words, a directed graph is rooted if it contains a directed spanning tree. 
An $n$-vertex directed tree is a rooted graph with $n-1$ arcs. It is easy to see that a directed tree has a unique root with an in-degree of $0$, while all other vertices have an in-degree of exactly $1$.
The smallest possible directed tree is a single isolated~vertex.

We use $a\bmod b$ to denote the modulo operation of two integers $a$ and $b$, which returns the remainder after dividing $a$ by $b$. 

\vspace{.05in}

{\em Algorithm 1:} Given $n>1$ vertices, label them, without loss of generality, from $1$ to $n$. Let $\bbb G(n,m)$ denote the $n$-vertex simple directed graph to be constructed with $m$ arcs. 
Start with the $m=n-1$ case and set $\bbb G(n,n-1)$ as any directed tree such that all its arcs $(i,j)$ satisfy $i<j$. 
For each integer $n\le m\le n(n-1)$, compute $v_{n,m} = n-((m-1) \bmod n)$, 
identify the smallest index $u_{n,m}\in\{1,\ldots,n\}$ such that $u_{n,m}\neq v_{n,m}$ and $(u_{n,m},v_{n,m})$ is not an arc in $\bbb G(n,m-1)$, 
then construct $\bbb G(n,m)$ by adding the arc $(u_{n,m},v_{n,m})$ to $\bbb G(n,m-1)$.
\hfill$\Box$

\vspace{.05in}

The requirement $i<j$ for all arcs $(i,j)$ immediately implies that each directed tree $\bbb G(n,n-1)$ is rooted at vertex~$1$ and has vertex $n$ as a leaf. 
Before proceeding, we first show that the above algorithm is well defined. 


\vspace{.05in}

\begin{lemma}\label{lm:available}
For any fixed $n>1$ and $n \le m \le n(n-1)$, the index $u_{n,m}$ defined in Algorithm 1 always exists.
\end{lemma}

\vspace{.05in}

{\bf Proof of Lemma \ref{lm:available}:}
From the algorithm description, $\bbb G(n,m)$ is constructed from $\bbb G(n,m-1)$ by adding an arc with head index $v_{n,m} = n-((m-1) \bmod n)$. This process proceeds inductively from $m=n$ to $m=n(n-1)$.
Note that for $n\le m\le n(n-1)$, the vertex index $v_{n,m} = n-((m-1) \bmod n)$ can take any value in $\{1,2,\ldots,n\}$. Each vertex in a graph with $n$ vertices has an in-degree of at most $n-1$.  
The algorithm begins with a directed tree, $\bbb G(n,n-1)$, in which vertex $1$ has in-degree $0$, while all other vertices have an in-degree of $1$. 
Hence, to prove the existence of the vertex index $u_{n,m}$ described in the algorithm, it is sufficient to show that 
vertex $1$ appears as the head of an added arc (i.e., $v_{n,m}=1$ occurs) at most $n-1$ times during the inductive construction process, while each vertex $i\in\{2,\ldots,n\}$ appears as the head of an added arc (i.e., $v_{n,m}=i$ occurs) at most $n-2$ times.
We thus consider vertex $1$ and vertices in $\{2,\ldots,n\}$ separately.  

First consider vertex 1. The condition $v_{n,m} = n-((m-1) \bmod n) = 1$ holds if and only if $m=kn$ with $k$ being any integer. Since $n \le m \le n(n-1)$, the condition occurs exactly $n-1$ times with $k\in\{1,\ldots,n-1\}$. 
Next consider vertices other than vertex 1. For each $i\in\{2,\ldots,n\}$, the condition $v_{n,m} = n-((m-1) \bmod n) = i$ holds if and only if $m=kn+1-i$ with $k$ being any integer. Since $n \le m \le n(n-1)$, the condition occurs exactly $n-2$ times with $k\in\{2,\ldots,n-1\}$. 
\hfill $\qed$

\vspace{.05in}

\begin{figure}[!ht]
\centering
\includegraphics[width=3.4in]{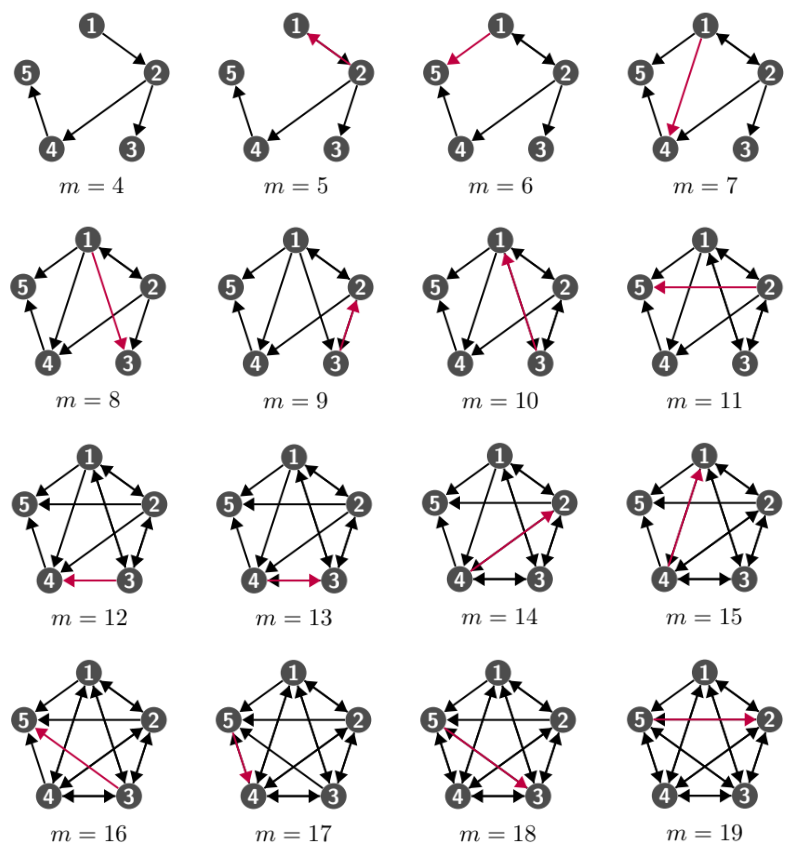} 
\caption{An example of the inductive construction process of Algorithm~1 for $n=5$ and $4\le m\le 19$}
\label{fig:algorithm5nodes}
\end{figure}

\begin{figure}[!ht]
\centering
\includegraphics[width=3.4in]{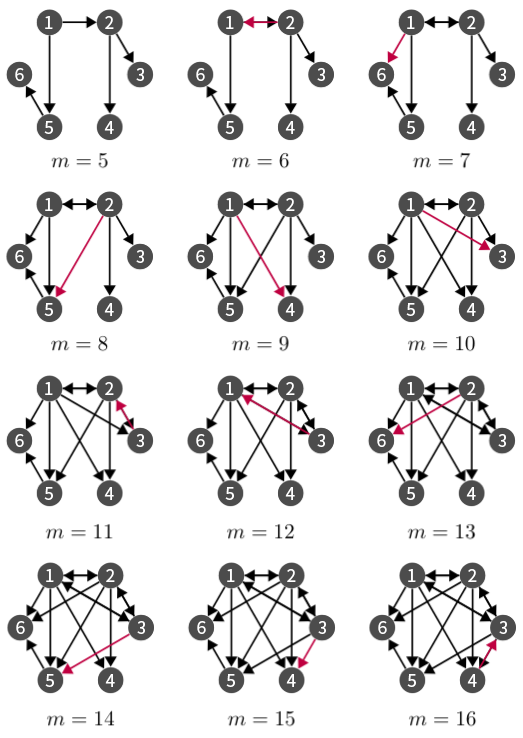} 
\caption{An example of the inductive construction process of Algorithm~1 for $n=6$ and $5\le m\le 16$}
\label{fig:algorithm6nodes}
\end{figure}

With the fact $v_{n,m} \in\{1,\ldots,n\}$, Lemma \ref{lm:available} ensures that Algorithm 1 operates without ambiguity under the given conditions for $n$ and $m$. 
In Figure~\ref{fig:algorithm5nodes}, we present an illustrative example of the algorithm that inductively constructs a sequence of graphs $\bbb G(n,m)$, where $n=5$ and $m$ ranges from $4$ to $19$. The first graph is a directed tree, and in each subsequent graph, a new arc is added to the preceding graph, with the newly added arc highlighted in purple.
We will consistently use purple to indicate newly added arcs when illustrating an inductive construction process for Algorithm 1. 
For simplicity in drawing, we use a bidirectional edge to represent two arcs in opposite directions throughout this paper; each bidirectional edge is therefore counted as two~arcs.
Another illustrative example for 6 vertices is provided in Figure \ref{fig:algorithm6nodes}.

To construct a graph with $n$ vertices and $m \ge n$ arcs using Algorithm~1, the computational complexity is $O(m-n+1)$.  This is because the algorithm initializes the graph as a directed tree with $n-1$ arcs and then builds it incrementally by adding the remaining $m-n+1$ arcs one at a time, with constant-time endpoint identification at each step.

\begin{figure}[!ht]
\centering
\includegraphics[width=3.45in]{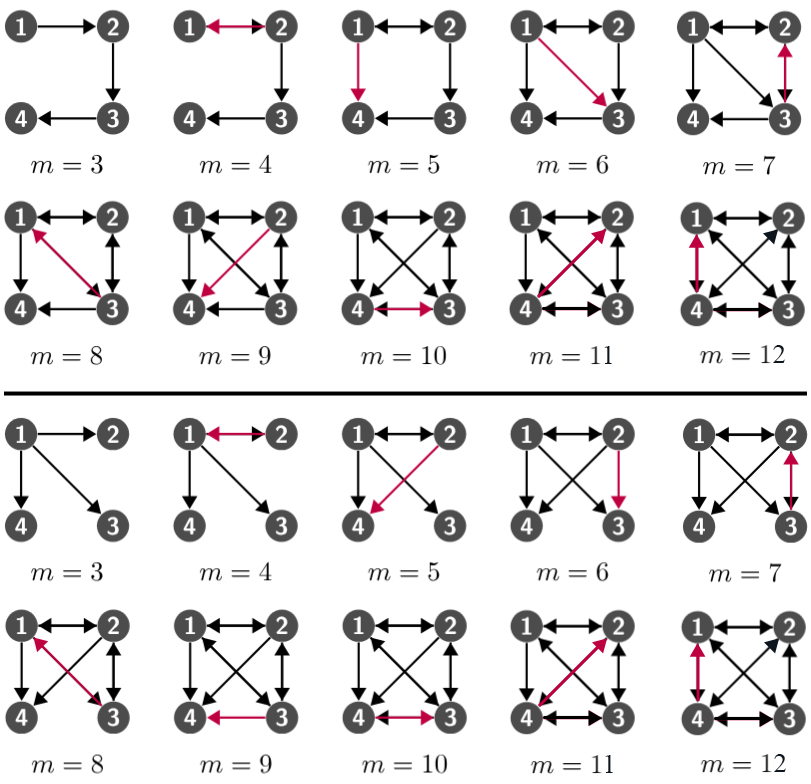} 
\caption{Two complete examples of the inductive construction process of Algorithm~1 for $n=4$}
\label{fig:algorithm4nodes2sequences}
\end{figure}

From the algorithm description, the constructed graphs are all rooted simple directed graphs. Indeed, any simple directed graph possessing the Laplacian spectrum specified in \eqref{eq:spectrum}, if it exists, is necessarily rooted (cf. Subsection \ref{subsec:rooted}). 
Moreover, each constructed graph $\bbb G(n,m)$, with $n\le m\le n(n-1)$, is dependent on the specific directed tree $\bbb G(n,n-1)$. In other words, for any fixed $n$, each distinct directed tree $\mathbb{G}(n, n-1)$ uniquely determines a corresponding sequence of graphs $\mathbb{G}(n, m)$, $n \leq m \leq n(n-1)$, with $\mathbb{G}(n,n(n-1))$ always being the complete graph. Figure \ref{fig:algorithm4nodes2sequences} illustrates two such complete sequences of graphs generated by the algorithm for $n=4$, starting from two different directed trees and ending at the same complete graph.

We say that two directed graphs are identical if they have the same sets of vertices and arcs, with the same labels on both; this requirement is stronger than graph isomorphism.
The classic Cayley's formula \cite{cayley} states that the number of distinct (undirected) trees on $n$ labeled vertices is $n^{n-2}$. 
Based on this, counting the number of distinct directed trees is easy. Any directed tree can be formed by orienting the edges of an undirected tree. Since each undirected tree allows any of its vertices to act as the root of a directed tree, each undirected tree can be oriented in $n$ different ways to form a directed tree. Therefore, the number of distinct directed trees on $n$ labeled vertices is $n\times n^{n-2}=n^{n-1}$. But this number cannot be used to count the total number of possible directed trees $\bbb G(n,n-1)$, 
as the algorithm requires that all arcs $(i,j)$ in $\bbb G(n,n-1)$ satisfy $i<j$.

\vspace{.05in}

\begin{lemma}\label{lm:number-of-directed-trees}
The number of distinct directed trees on $n$ labeled vertices, such that each arc $(i,j)$ satisfies $i<j$, is $(n-1)!$.
\end{lemma}

\vspace{.05in}

{\bf Proof of Lemma \ref{lm:number-of-directed-trees}:}
We prove the lemma by induction on $n$. For the base case $n=1$, the lemma is clearly true. For the inductive step, suppose that the lemma holds for $n=k$, where $k$ is a positive integer. Let $n=k+1$. Since each arc $(i,j)$ is required to satisfy $i<j$, vertex $k+1$ must be a leaf vertex, and its parent vertex can be any vertex in $\{1,\ldots,k\}$. Thus, for each directed tree with $k$ vertices that satisfies the requirement, there are $k$ different ways to construct a directed tree with $k+1$ vertices that also satisfies the requirement by adding vertex $k+1$ as a child of any existing vertex. By the inductive hypothesis, the total number of desired directed trees with $k$ vertices is $(k-1)!$. Therefore, the total number of such directed trees with $k+1$ vertices equals $k\times (k-1)! = k!$, which proves the inductive step. 
\hfill $\qed$

\vspace{.05in}

\begin{figure}[!ht]
\centering
\includegraphics[width=3.4in]{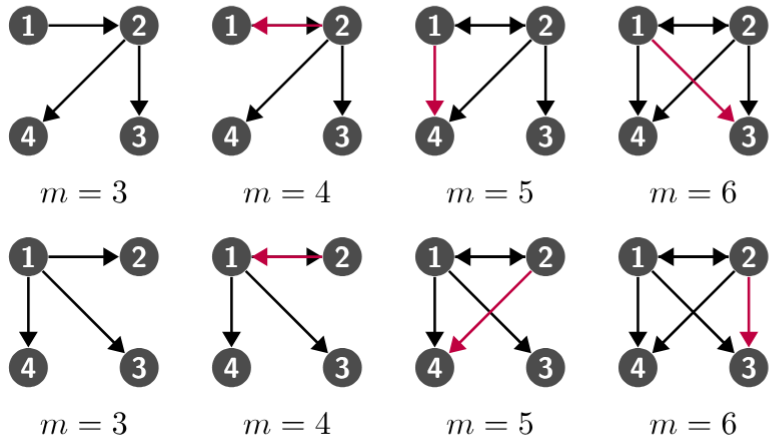} 
\caption{Two inductive construction processes (first row and second row) of Algorithm~1 for $n=4$, starting from different directed trees and leading to the identical $\bbb G(4,6)$ graph}
\label{fig:algorithm-identical-graphs}
\end{figure}

The lemma states that there are $(n-1)!$ different possible $\mathbb{G}(n,n-1)$. 
That is to say, each $\mathbb{G}(n,m)$ may represent up to $(n-1)!$ different graphs. 
For simplicity, we use the notation as is and take this fact without further mention in the sequel. 
It is possible that, for certain values of $n$ and $m$, the graph $\bbb G(n,m)$ constructed by the algorithm may be identical, even when the construction process begins with different directed trees. 
A trivial example is when $m=n(n-1)$ with which $\bbb G(n,m)$ must be the complete graph regardless of $\bbb G(n,n-1)$. Another illustrative $4$-vertex example is given in Figure~\ref{fig:algorithm-identical-graphs}.
Figure \ref{fig:algorithm4nodes2sequences} provides a further $4$-vertex example, where two construction processes start from two different directed trees but end at the same graph $\mathbb{G}(4,m)$ for $m \in \{9,10,11,12\}$. This example highlights that once two construction processes for the same vertex number $n$ coincide at some $\mathbb{G}(n,m)$, all subsequent graphs in the two sequences must also coincide. This follows directly from the algorithm description.
In addition, 
a sufficient condition on the relationship between $n$ and $m$ is provided in Lemma~\ref{prop:G-large-m}, which guarantees that 
$\bbb G(n,m)$ is unique no matter what the initial directed tree is.


All constructed $\bbb G(n,m)$ graphs have the following property. 
A directed graph is called almost regular if the difference between its largest and smallest in-degrees is at most $1$.
In the special case when all in-degrees are equal, the graph is called regular.
Let $\nu \dfb \lfloor \frac{m}{n} \rfloor$.

\vspace{.05in}

\begin{proposition}\label{prop:Galmostregular}
For any integers $n\ge 2$ and $m\ge n-1$, any $\bbb G(n,m)$ constructed by Algorithm 1 is almost regular, with $n(\nu+1)-m$ vertices of in-degree $\nu$ and $m-n\nu$ vertices of in-degree $\nu+1$; 
that is, its vertex in-degree sequence is 
\eq{(d_1, \ldots, d_n) = (\underbrace{\nu,\ldots,\nu}_{n(\nu+1)-m}, \; \underbrace{\nu+1,\ldots,\nu+1}_{m-n\nu}\;).\label{eq:indegree}} 
\end{proposition}

\vspace{.1in}

The proposition implies that whenever $m$ is a multiple of $n$, any graph $\bbb G(n,m)$ constructed by Algorithm 1 is an exactly regular directed graph.

\vspace{.05in}

{\bf Proof of Proposition \ref{prop:Galmostregular}:}
We first consider two special cases. 
First, in the case when $m = n-1$, $\bbb G(n,n-1)$ is a directed tree, which is clearly almost regular. 
Second, in the case when $m=n$, $\bbb G(n,n)$ is constructed by adding the arc $(u_{n,n}, 1)$ to $\bbb G(n,n-1)$. Then, all $n$ vertices have an in-degree of $1$, and thus $\bbb G(n,n)$ is regular. 
It remains to consider the case $m \geq n+1$.
From the algorithm description, $\bbb G(n,m)$ is constructed from $\bbb G(n,n)$ by inductively adding arcs whose head indices are given by $v_{n,k} = n-((k-1) \bmod n)$ for $k\in\{n+1,\ldots,m\}$. Note that as $k$ ranges from $n+1$ to $2n$, $v_{n,k}$ takes values from $n$ to $1$, and this pattern repeats with a period of $n$ as $k$ continues from $2n+1$ to $m$. This implies that the in-degrees of $\bbb G(n,m)$ satisfy $d_n \ge \cdots \ge d_1$ and $d_n-d_1\le 1$ for all $m\ge n+1$. Thus, $\bbb G(n,m)$ is always almost regular. 
The remaining statement of the proposition directly follows from the following lemma.~\hfill $\qed$

\vspace{.05in}

\begin{lemma}\label{lem:degreesequence}
    For any almost regular simple directed graph with $n$ vertices and $m$ arcs, assume, without loss of generality, that its vertex in-degrees satisfy $d_1\le \cdots \le d_n$. Then, its in-degree sequence is \eqref{eq:indegree}.
\end{lemma}

\vspace{.05in}

{\bf Proof of Lemma \ref{lem:degreesequence}:}
Since the graph is almost regular, $d_n-d_1\le 1$. 
Suppose there are $1\le p\le n$ vertices with the minimal in-degree $d_1$. Then, the remaining $q=n-p$ vertices have an in-degree of $d_1+1$. It follows that $m=pd_1+q(d_1+1)=d_1n+q$. As $q$ takes a value in $\{0,1,\ldots,n-1\}$, $d_1$ and $q$ are respectively the unique quotient and remainder when $m$ is divided by $n$. Then, $d_1=\lfloor \frac{m}{n} \rfloor = \nu$ and $q=m-d_1n=m-n\nu$. Therefore, the in-degree sequence is \eqref{eq:indegree}.
\hfill $\qed$

\vspace{.05in}

Proposition \ref{prop:Galmostregular} specifies the in-degree $d_i$ of each vertex $i$ in $\bbb G(n,m)$. Each $d_i$ can take an integer value from 0 to $n-1$, depending on the value of $m$. The following lemma further identifies the in-neighbors corresponding to these $d_i$ values.

\vspace{.05in}

\begin{lemma}\label{lem:incoming-neighbors}
    For any $\bbb G(n,m)$, the $d_1$ incoming arcs of vertex~$1$ originate from $d_1$ vertices whose indices are in $\{2,\ldots,d_1+1\}$, and for each $i\in\{2,\ldots,n\}$, the $d_i$ incoming arcs of vertex $i$ originate from vertex $i_j$, the unique vertex index such that $(i_j,i)$ is an arc in $\bbb G(n,n-1)$, and from $d_i-1$ other vertices whose indices are the $d_i-1$ smallest elements of $\{1,\ldots,n\} \setminus \{i, i_j\}$.
\end{lemma}

\vspace{.05in}

{\bf Proof of Lemma \ref{lem:incoming-neighbors}:}
From the algorithm description, vertex $1$ has in-degree $0$ in $\bbb G(n,n-1)$, so each of its incoming arcs in $\bbb G(n,m)$, if any, must be added as $(u_{n,p},v_{n,p})$ with $v_{n,p}=1$ for some index $p\in\{n,\ldots,m\}$ during the inductive construction process. 
Since each $u_{n,p}$ is defined as the smallest vertex index such that $u_{n,p} \neq 1$ and $(u_{n,p},1)$ is not an arc in $\bbb G(n, p-1)$, the $d_1$ incoming neighbors of vertex $1$ must be the $d_1$ vertices with the smallest indices other than $1$, that is, the vertices in $\{2,\ldots,d_1+1\}$. 

Next consider any vertex $i\in\{2,\ldots,n\}$, which has exactly one incoming arc $(i_j,i)$ in the directed tree $\bbb G(n,n-1)$, and thus has in-degree at least one in $\bbb G(n,m)$.
The remaining incoming arcs of vertex $i$ in $\bbb G(n,m)$, if any, are added through the inductive construction process. Using the same argument as in the previous paragraph, the $d_i-1$ remaining incoming arcs originate from the $d_i-1$ vertices with the smallest indices in $\{1,\ldots,n\} \setminus \{i, i_j\}$.    
\hfill$\qed$

\vspace{.05in}

The most important property of $\bbb G(n,m)$ is stated below. Additional properties will be presented later in the paper.

\vspace{.05in}

\begin{theorem}\label{th:pnas-spectrum}
For any integers $n\ge 2$ and $m\ge n-1$, any graph constructed by Algorithm 1, $\bbb G(n,m)$, has the Laplacian spectrum given in \eqref{eq:spectrum}. 
\end{theorem}

\vspace{.05in}


To prove the theorem, we need several concepts and results. 

A directed graph is acyclic if it contains no directed cycles. Thus, by definition, a directed acyclic graph cannot contain a self-arc. Any directed tree is acyclic. The transpose of an acyclic graph remains acyclic.

\vspace{0.05in}

\begin{lemma}\label{lm:acyclic}
For any acyclic simple directed graph, its Laplacian spectrum consists of its in-degrees.
\end{lemma}

\vspace{.05in}

{\bf Proof of Lemma \ref{lm:acyclic}:}
The adjacency matrix $A$ of a directed graph $\bbb G$, as defined in the introduction, is based on in-degrees. The out-degree based adjacency matrix is the transpose of the in-degree based adjacency matrix; that is, its $ij$th entry equals 1 if $(i,j)$ is an arc in the graph, and equals 0 otherwise, as referenced in \cite[Page 151]{Harary69}. 
For any permutation matrix $P$, $P'AP$ represents an adjacency matrix of the same graph, but with its vertices relabeled; the same property applies to out-degree based adjacency matrices.
Since $\bbb G$ is acyclic, from \cite[Theorem 16.3]{Harary69}, there exists a permutation matrix $P$ with which $P'A'P$ is upper triangular. Then, $P'AP$ is lower triangular, which implies that $P'LP$ is also lower triangular. Thus, the spectrum of $P'LP$ consists of its diagonal entries. Since $P'LP$ and $L$ share the same spectrum and diagonal entries, the spectrum of $L$ consists of its diagonal entries, which are the in-degrees of $\bbb G$. 
\hfill $\qed$

\vspace{.05in}

The union of two directed graphs, $\bbb G_1$ and $\bbb G_2$, with the same vertex set, denoted by $\bbb G_1\cup\bbb G_2$, is the directed graph with the same vertex set and its directed edge set being the union of the directed edge sets of $\bbb G_1$ and $\bbb G_2$. Similarly, the intersection of two directed graphs, $\bbb G_1$ and $\bbb G_2$, with the same vertex set, denoted by $\bbb G_1\cap\bbb G_2$, is the directed graph with the same vertex set and its directed edge set being the intersection of the directed edge sets of $\bbb G_1$ and $\bbb G_2$.
Graph union is an associative binary operation, and thus the definition extends unambiguously to any finite sequence of directed graphs.
The complement of a simple directed graph $\bbb G$, denoted by $\overline{\bbb G}$, is the simple directed graph with the same vertex set such that $\bbb G \cup \overline{\bbb G}$ equals the complete graph and $\bbb G \cap \overline{\bbb G}$ equals the empty graph. 
It is easy to see that if vertex $i$ has in-degree $d_i$ in $\bbb G$, then it has in-degree $n-1-d_i$ in $\overline{\bbb G}$. Moreover, the total number of arcs in $\bbb G$ and $\overline{\bbb G}$ is $n(n-1)$.



It is easy to see that any Laplacian matrix has an eigenvalue at $0$ with an eigenvector $\1$, where $\1$ is the column vector in $\R^n$ with all entries equal to $1$. Using the Gershgorin circle theorem \cite{gershgorin}, it is straightforward to show that all Laplacian eigenvalues, except for those at 0, have positive real parts, as was done in \cite[Theorem 2]{reza1} for out-degree based Laplacian matrices. Thus, the smallest real part of any Laplacian eigenvalue is always $0$. 
More can be said.

\vspace{.05in}

\begin{lemma}\label{lm:complement}
If the Laplacian spectrum of an $n$-vertex simple directed graph $\bbb G$ is $\{0,\lambda_2, \ldots, \lambda_n\}$ with $0\le {\rm Re}(\lambda_2)\le \cdots \le {\rm Re}(\lambda_n)$, then the Laplacian spectrum of its complement $\overline{\bbb G}$ is $\{0, n-\lambda_n, \ldots, n-\lambda_2\}$ and $0\le {\rm Re}(n-\lambda_n)\le \cdots \le {\rm Re}(n-\lambda_2)$.
\end{lemma}

\vspace{.05in}

The following proof of the lemma employs the same technique as that used in the proof of Theorem 2 in \cite{Agaev05}, which was developed for a variant of Laplacian matrices.
For any square matrix $M$, we denote its characteristic polynomial as $p_M(\lambda) = \det(\lambda I - M)$ in the sequel.

\vspace{.05in}

{\bf Proof of Lemma \ref{lm:complement}:}
Let $L$ and $\overline{L}$ be the Laplacian matrices of $\bbb G$ and $\overline{\bbb G}$, respectively. It is straightforward to verify that $L+\overline{L} = nI-J$, where $I$ is the identity matrix and $J$ is the $n \times n$ matrix with all entries equal to 1.
Let $Q = L+J = nI-\overline{L}$. 
We first show that 
\begin{equation}\label{eq:complement-pf}
\lambda p_Q(\lambda) = (\lambda-n) p_L(\lambda).
\end{equation}
Note that $p_L(\lambda)=0$ when $\lambda=0$, as $L$ has an eigenvalue at 0. Thus, \eqref{eq:complement-pf} holds when $\lambda=0$.

To prove \eqref{eq:complement-pf} for $\lambda\neq 0$, let $c_i$, $i\in\{1,\ldots,n\}$ denote the $i$th column of matrix $\lambda I - L$. Since $Q = L+J$, it follows that the $i$th column of matrix $\lambda I-Q$ is $c_i-\1$. 
Since the determinant of a matrix is multilinear and adding one column to another does not alter its value, 
$p_Q(\lambda) = \det\; [c_1-\1, c_2-\1, \cdots, c_n-\1] = \det\; [c_1, c_2-\1, \cdots, c_n-\1] - \det\; [\1, c_2-\1, \cdots, c_n-\1] = \det\; [c_1, c_2-\1, \cdots, c_n-\1] - \det\; [\1, c_2, \cdots, c_n]$. Repeating this process sequentially for the columns from 2 to $n$ leads to
\begin{align*}
    p_Q(\lambda) = p_L(\lambda) 
    - \sum_{i=1}^n \det \big[c_1, \cdots, c_{i-1},\1, c_{i+1},\cdots,c_n\big].
\end{align*}
Note that $\sum_{j=1}^n c_j =\lambda\1$, as each row sum of $\lambda I - L$ is equal to $\lambda$. Then, for any $i\in\{1,\ldots,n\}$,
\begin{align*}
    &\det \big[c_1, \cdots, c_{i-1},\1, c_{i+1},\cdots,c_n\big] \\
    =\;& \det \big[c_1, \cdots, c_{i-1},\textstyle\frac{1}{\lambda}\sum_{j=1}^n c_j, c_{i+1},\cdots,c_n\big] \\
    =\;& \textstyle\frac{1}{\lambda}\det \big[c_1, \cdots, c_{i-1},\sum_{j=1}^n c_j, c_{i+1},\cdots,c_n\big] \\
    =\;& \textstyle\frac{1}{\lambda}\det \big[c_1, \cdots, c_n\big] = \frac{1}{\lambda} p_L(\lambda).
\end{align*}
Substituting this equality into the preceding expression for $p_Q(\lambda)$ yields $p_Q(\lambda) = p_L(\lambda) - \frac{n}{\lambda}p_L(\lambda)$, which proves \eqref{eq:complement-pf}.

To proceed, recall that $\overline L = nI-Q$. Then, 
\begin{align*}
p_{\overline{L}}(\lambda) &= \det(\lambda I - \overline{L}) = \det(\lambda I - nI + Q) \\
&= (-1)^n \det((n-\lambda)I - Q) = (-1)^n p_Q(n-\lambda). 
\end{align*}
From this and \eqref{eq:complement-pf}, with $\lambda$ substituted by $n-\lambda$,
\begin{align}\label{eq:eigrelation}
    (n-\lambda)p_{\overline{L}}(\lambda) = (-1)^{n+1} \lambda p_L(n-\lambda).
\end{align}
Both sides of \eqref{eq:eigrelation} are polynomials in $\lambda$ of degree $n+1$. It is easy to see that 0 and $n$ are roots of both sides, as 0 is an eigenvalue of both $L$ and $\overline{L}$. Then, the set of nonzero roots of $p_{\overline{L}}(\lambda)$ coincides with the set of roots of $p_L(n-\lambda)$, excluding the root at $n$.
Therefore, if the Laplacian spectrum of $\bbb G$ is $\{0,\lambda_2, \ldots, \lambda_n\}$, then the Laplacian spectrum of $\overline{\bbb G}$ is $\{0, n-\lambda_n, \ldots, n-\lambda_2\}$. Recall that the smallest real part of all Laplacian eigenvalues is always 0. With these facts, it is easy to see that if $0\le {\rm Re}(\lambda_2)\le \cdots \le {\rm Re}(\lambda_n)$, then $0\le {\rm Re}(n-\lambda_n)\le \cdots \le {\rm Re}(n-\lambda_2)$.
\hfill $\qed$

\vspace{.05in}

The disjoint union of two directed graphs is a larger directed graph whose vertex set is the disjoint union of their vertex sets, and whose arc set is the disjoint union of their arc sets. Disjoint union is an associative binary operation, and thus the definition extends unambiguously to any finite sequence of directed graphs. Any disjoint union of two or more graphs is necessarily disconnected.
A directed forest is a disjoint union of directed tree(s). A directed forest composed of $k$ directed trees thus has $k$ vertices with an in-degree of 0, while all other vertices have an in-degree of exactly 1.
It is easy to see that the number of directed trees in a directed forest is equal to the difference between the number of vertices and the number of arcs. 
An $n$-vertex directed star is a directed tree whose $n-1$ arcs all originate from the root.

\vspace{.05in}

\begin{lemma}\label{prop:G-large-m}
For any integers $n$ and $m$ such that $n \ge 2$ and $(n-1)^2 \le m \le n(n-1)$, any $\bbb G(n,m)$ constructed by Algorithm 1 is the complement of the $n$-vertex directed forest consisting of a directed star with $(n(n-1)-m+1)$ vertices, rooted at vertex $n$ and with leaf vertices $1$ to $(n(n-1)-m)$, and $m-(n-1)^2$ isolated vertices.
\end{lemma}

\vspace{.05in}


\begin{figure}[!ht]
\centering
\includegraphics[width=1.8in]{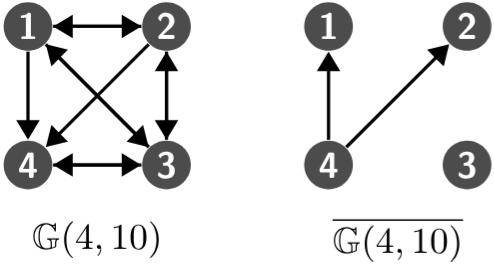} 
\caption{Example illustrating Lemma \ref{prop:G-large-m}}
\label{fig:prop-G-large-m}
\end{figure}

The condition $m\ge (n-1)^2$ ensures $n(n-1)-m\le n$. In the special case when $m=n(n-1)$, the directed star reduces to an isolated vertex, which is consistent with the fact that the complement of a complete graph is an empty graph. 
Figure \ref{fig:prop-G-large-m} provides a simple example to illustrate Lemma \ref{prop:G-large-m}. Note that $\bbb G(4,10)$, as shown in the figure, is unique. Indeed, 
Lemma~\ref{prop:G-large-m} immediately implies that if $(n-1)^2 \le m \le n(n-1)$, then $\bbb G(n,m)$ is unique regardless of which directed tree $\bbb G(n,n-1)$ is. 
Figure \ref{fig:algorithm4nodes2sequences} provides two complete sequences of $\bbb G(n,m)$ graphs for $n=4$, inductively constructed by Algorithm 1. The two sequences start from different directed trees, and from $m=9$ onward, the constructed graphs become identical. 
Note that the uniqueness of $\bbb G(4,9)$, $\bbb G(4,10)$, and $\bbb G(4,11)$ can also be implied by Proposition \ref{prop:Galmostregular} and Lemma \ref{lem:incoming-neighbors}.


To prove Lemma \ref{prop:G-large-m}, we need the following result. 


\vspace{.05in}

\begin{lemma}\label{lm:deg-n-1}
If there exists an index $k\in\{1,\ldots,n-1\}$ such that $(n,k)$ is an arc in $\bbb G(n,m)$, then the in-degree of vertex $k$ in $\bbb G(n,m)$ is $n-1$. 
\end{lemma}

\vspace{.05in}

{\bf Proof of Lemma \ref{lm:deg-n-1}:}
Note that $(n,k)$ cannot be an arc in $\bbb G(n,n-1)$ because any arc in $\bbb G(n,n-1)$ is required to have a head index greater than the tail index.
Since the arc $(n,k)$ is not in $\bbb G(n,n-1)$ but is in $\bbb G(n,m)$, from the algorithm description, there must exist exactly one index $p\in\{n,\ldots,m\}$ such that $u_{n,p}=n$ and $v_{n,p}=k$. The algorithm sets $u_{n,p}$ as the smallest vertex index for which $u_{n,p}\neq k$ and $(u_{n,p},k)$ is not an arc in $\bbb G(n,p-1)$. Thus, $(i,k)$ is an arc in $\bbb G(n,p-1)$ for all $i\in\{1,\ldots,n-1\}\setminus\{k\}$. This implies that the in-degree of vertex $k$ in $\bbb G(n,p-1)$ is $n-1$, and the same holds for $\bbb G(n,m)$. 
\hfill $\qed$

\vspace{.05in}

{\bf Proof of Lemma \ref{prop:G-large-m}:}
In the special case when $m=n(n-1)$, $\bbb G(n,m)$ is the complete graph, which is the complement of the empty graph. It is easy to verify that the lemma is true in this case. 
Next consider the case when $(n-1)^2 \le m < n(n-1)$ for which $\nu = \lfloor \frac{m}{n} \rfloor = n-2$. 
From Proposition \ref{prop:Galmostregular}, the vertex in-degree sequence of $\bbb G(n,m)$ is
$$(d_1, \ldots, d_n) = (\; \underbrace{n-2, \ldots, n-2}_{n(n-1)-m},\; \underbrace{n-1, \ldots, n-1}_{m-n(n-2)} \;).$$
Since $m\ge (n-1)^2$, $m-n(n-2)\ge 1$, which implies that the in-degree of vertex $n$ is $n-1$. 
Consequently, in $\overline{\bbb G(n,m)}$, the complement of $\bbb G(n,m)$, vertices from $1$ to $(n(n-1)-m)$ have an in-degree of $1$, and the remaining vertices, including vertex $n$, have an in-degree of $0$. 
Note that the total number of arcs in $\overline{\bbb G(n,m)}$ is $(n(n-1)-m)$. 
We claim that these $(n(n-1)-m)$ arcs form a directed star rooted at vertex $n$, with leaf vertices labeled from $1$ to $(n(n-1)-m)$, which immediately implies the lemma.  To prove the claim, it is equivalent to show that $(n,k)$ is an arc in $\overline{\bbb G(n,m)}$ for all $k\in\{1,\ldots,n(n-1)-m\}$. 
To this end, suppose to the contrary that $(n,k)$ is not an arc in $\overline{\bbb G(n,m)}$ for some $k\in\{1,\ldots,n(n-1)-m\}$. This implies that $(n,k)$ is an arc in $\bbb G(n,m)$. From Lemma \ref{lm:deg-n-1}, the in-degree of vertex $k$ in $\bbb G(n,m)$ is $n-1$.  
But this contradicts the fact that $d_k=n-2$. Therefore, the claim is true. 
\hfill $\qed$

\vspace{.05in}

Let $L(n,m)$ denote the Laplacian matrix of $\bbb G(n,m)$, with its $ij$th entry denoted by $[L(n,m)]_{ij}$. It has the following entry-wise property.

\vspace{.05in}

\begin{lemma}\label{lem:Lentries}
For any integers $n$ and $m$ such that $n \ge 2$ and $n-1 \le m \le (n-1)^2$, $[L(n,m)]_{in}=0$ for each $i\in\{1,\ldots,n-1\}$.
\end{lemma}

\vspace{.05in}

{\bf Proof of Lemma \ref{lem:Lentries}:}
From the definition of a Laplacian matrix, to prove the lemma, it is equivalent to showing that, for each $i \in \{1, \ldots, n-1\}$, the arc $(n, i)$ does not belong to $\bbb G(n, m)$ when $n-1 \le m \le (n-1)^2$.
From the inductive construction described in Algorithm~1, if arc $(n,i)$ belongs to a graph $\mathbb{G}(n,m)$ for some $n-1 \le m < (n-1)^2$, then there must exist a graph $\mathbb{G}(n,(n-1)^2)$ that also contains arc $(n,i)$. Therefore, it is sufficient to show that any graph $\mathbb{G}(n,(n-1)^2)$ constructed by Algorithm~1 does not contain the arc $(n,i)$.
From Lemma \ref{prop:G-large-m}, any $\mathbb{G}(n,(n-1)^2)$ is the complement of the $n$-vertex directed star rooted at vertex~$n$.
It follows that $\mathbb{G}(n,(n-1)^2)$ contains no arcs of the form $(n,i)$ for any $i \in \{1,\ldots,n-1\}$.
\hfill$\qed$

\vspace{.05in}

Lemma \ref{lem:Lentries} has the following important implication. 
Recall that Algorithm 1 may generate different graph sequences depending on the initial tree graph. 

\vspace{.05in}

\begin{lemma}\label{lm:delete-vertex-n}
Let $\bbb G(n,m)$ be the graph constructed by Algorithm 1 starting from a directed tree $\bbb G(n,n-1)$. 
If $n \ge 3$ and $n-1 \le m \le (n-1)^2$, then the subgraph of $\bbb G(n,m)$ induced by the vertex subset $\{1,\ldots,n-1\}$ is $\bbb G(n-1, m-d_n)$, the graph constructed by Algorithm 1 starting from the subgraph of $\bbb G(n,n-1)$ induced by the vertex subset $\{1,\ldots,n-1\}$, where $d_n$ is the in-degree of vertex $n$ in $\bbb G(n,m)$.  
\end{lemma}


\vspace{.05in}


Since $\bbb G(n,n-1)$ is an $n$-vertex directed tree with vertex $n$ as a leaf, its subgraph induced by the vertex subset $\{1,\ldots,n-1\}$ is a directed tree with $n-1$ vertices. 
The subgraph also satisfies the requirement that every arc has a head index greater than the tail index, so it can serve as an initial tree graph $\bbb G(n-1,n-2)$ for Algorithm 1 to construct graphs with $n-1$ vertices.
In addition, recall that $\bbb G(n,n-1)$ is rooted at vertex $1$, and so is $\bbb G(n,m)$ for any $m\in\{n,\ldots,n(n-1)\}$. Then, the subgraph of $\bbb G(n,m)$ induced by the vertex subset $\{1,\ldots,n-1\}$ is also rooted at vertex $1$ and therefore has at least $n-2$ arcs. Since $m-d_n$ is an upper bound on the total number of arcs in the subgraph, its value is no smaller than $n-2$. In the following proof, we will soon show that $m-d_n\le (n-1)(n-2)$. These two facts guarantee that $\bbb G(n-1,m-d_n)$ is well-defined.



\vspace{.05in}

\textbf{Proof of Lemma \ref{lm:delete-vertex-n}:} Let $\bbb H$ be the subgraph of $\bbb G(n,m)$ induced by the vertex subset $\{1,\ldots,n-1\}$. From Lemma~\ref{lem:Lentries}, vertex $n$ has no outgoing arcs in $\bbb G(n,m)$. Then, $\bbb H$ has $n-1$ vertices and $m-d_n$ arcs, which implies $m-d_n\le (n-1)(n-2)$.  
In addition, each vertex $i\in\{1,\ldots,n-1\}$ has the same in-degree $d_i$ in $\bbb H$ as in $\bbb G(n,m)$, with the values of $d_i$, $i\in\{1,\ldots,n\}$, being given in \eqref{eq:indegree}. 
Let $b_i$ denote the in-degree of vertex $i$ in $\bbb G(n-1,m-d_n)$. 
From Proposition~\ref{prop:Galmostregular}, 
$$(b_1, \ldots, b_{n-1}) = (\underbrace{u, \ldots, u,}_{(n-1)(u+1)-(m-d_n)} \underbrace{u+1, \ldots, u+1}_{(m-d_n) - (n-1)u} \;\;), $$
where $u = \lfloor \frac{m-d_n}{n-1} \!\rfloor$. We claim that $b_i=d_i$ for all $i\in\{1,\ldots,n-1\}$. To prove the claim, we consider two scenarios separately.
First, suppose that $n$ divides $m$. Then, from \eqref{eq:indegree}, all $d_i$, $i\in\{1,\ldots,n\}$ equal $\nu=\frac{m}{n}$, which implies that $n-1$ divides $m-d_n$ and $u=\nu$. It follows that all $b_i$, $i\in\{1,\ldots,n-1\}$ equal $u$, and thus the claim holds. 
Next, suppose that $n$ does not divide $m$. Then, $d_n=\nu +1=\lfloor\frac{m}{n}\rfloor+1$ and $m=\nu n+r$, where $1\le r\le n-1$ is the unique remainder when $m$ is divided by $n$. With these, $u = \lfloor \frac{m-d_n}{n-1} \!\rfloor = \lfloor \frac{\nu n+r-\nu -1}{n-1}\rfloor = \nu + \lfloor \frac{r-1}{n-1} \rfloor = \nu$ and thus $(n-1)(u+1)-(m-d_n)=n(\nu+1)-m$, which validates the claim. 
This ensures that for each $i\in\{1,\ldots,n-1\}$, vertex~$i$ has an in-degree of $d_i$ in both $\bbb H$ and $\bbb G(n-1,m-d_n)$. 
Since $\bbb G(n-1,n-2)$ is the subgraph of $\bbb G(n,n-1)$ induced by the vertex subset $\{1,\ldots,n-1\}$, where the former is a directed tree with $n-1$ vertices and the latter is a directed tree with $n$ vertices, for each $i\in\{1,\ldots,n-1\}$, vertex $i$ has the same unique in-neighbor index in both directed trees. With this, Lemma \ref{lem:incoming-neighbors} implies that each vertex has the same set of in-neighbor indices in $\bbb H$ and $\bbb G(n-1,m-d_n)$. Therefore, $\bbb H = \bbb G(n-1,m-d_n)$. 
\hfill $\qed$

\vspace{0.05in}


We will also need the following lemmas regarding the relationship between $\kappa=\lfloor \frac{m}{n-1} \rfloor$ and $\nu=\lfloor \frac{m}{n} \rfloor$. 

\vspace{.05in}

\begin{lemma}\label{lm:p-k}
$\nu \in \{\kappa-1,\kappa\}$ for any integers $n$ and $m$ such that $n \ge 2$ and $1\le m \le n(n-1)$.
\end{lemma}

\vspace{.05in}

{\bf Proof of Lemma \ref{lm:p-k}:} Since $\lfloor \frac{m}{n-1} \rfloor \le \frac{m}{n-1}$ and $\lfloor \frac{m}{n} \rfloor > \frac{m}{n}-1$, it follows that 
$\kappa-\nu = \lfloor \frac{m}{n-1} \rfloor - \lfloor \frac{m}{n} \rfloor < \frac{m}{n-1} - (\frac{m}{n}-1 ) = \frac{m}{n(n-1)} + 1 \le 2$. 
As $\kappa-\nu$ is a nonnegative integer, it can only take a value of either 0 or 1, which implies that $\nu$ is equal to either $\kappa$ or $\kappa-1$.
\hfill $\qed$

\vspace{0.05in}

\begin{lemma}\label{lm:floor}
$\lfloor \frac{m-\nu-1}{n-2} \rfloor = \kappa$ for any integers $n$ and $m$ such that 
$n \ge 3$ and $1 \le m \le n(n-2)$.
\end{lemma}

\vspace{0.05in}

{\bf Proof of Lemma \ref{lm:floor}:}
First, consider the special case when $1 \le m \le n-2$, which implies $\kappa = \nu = 0$. Then, $\lfloor \frac{m-\nu-1}{n-2} \rfloor = \lfloor \frac{m-1}{n-2} \rfloor = 0 = \kappa$. Thus, the lemma holds in this case. 

Next, consider the general case when 
$n-1 \le m \le n(n-2)$.
Note that $m$ can be written as $m = \kappa(n-1)+r$, where $0 \le r \le n-2$ is the unique remainder when $m$ is divided by $n-1$. 
From Lemma \ref{lm:p-k}, $\nu$ equals either $\kappa$ or $\kappa-1$. Let us first suppose $\nu=\kappa$. Then, $\frac{\kappa (n-1) + r}{n} = \frac{m}{n} \ge \lfloor \frac{m}{n} \rfloor =\nu =\kappa$, which implies $r \ge \kappa = \lfloor \frac{m}{n-1} \rfloor \ge 1$. Thus, 
$\lfloor \frac{m-\nu-1}{n-2} \rfloor = \lfloor \frac{\kappa(n-1)+r-\kappa-1}{n-2} \rfloor = \kappa + \lfloor \frac{r-1}{n-2} \rfloor = \kappa$. 
In the next step, we suppose $\nu=\kappa-1$. Then, $\lfloor \frac{m-\nu-1}{n-2} \rfloor = \lfloor \frac{\kappa(n-1)+r-\kappa}{n-2} \rfloor = \kappa + \lfloor \frac{r}{n-2} \rfloor$, which equals $\kappa$ if $0\le r< n-2$. 
To complete the proof, it remains to consider the case when $r=n-2$. We claim that $r\neq n-2$. To prove the claim, suppose to the contrary that $r = n-2$, with which $m = \kappa(n-1)+r = \kappa n+(n-2-\kappa)$.
Meanwhile, 
as $n \ge 3$ and $1 \le m \le n(n-2)$,
$\kappa = \lfloor \frac{m}{n-1} \rfloor \le \frac{n(n-2)}{n-1}<n-1$, which implies that $n-1-\kappa$ is a positive integer. Then, 
$n-2-\kappa \ge 0$, and thus $\nu = \lfloor \frac{m}{n} \rfloor = \lfloor \frac{\kappa n+(n-2-\kappa)}{n} \rfloor = \kappa+\lfloor \frac{n-2-\kappa}{n} \rfloor = \kappa$. 
But this contradicts $\nu=\kappa-1$. Therefore, $r\neq n-2$. 
\hfill $\qed$


\begin{lemma}\label{lm:n-divides-m}
For any integers $n$ and $m$ such that 
$n \ge 3$ and $1 \le m < n(n-2)$,
if $n$ divides $m$, then $\kappa=\nu=\lfloor \frac{m-\nu}{n-2} \rfloor$.
\end{lemma}

\vspace{.05in}

{\bf Proof of Lemma \ref{lm:n-divides-m}:}
Since $n$ divides $m$, $m=\nu n$. 
Then, $\kappa = \lfloor \frac{m}{n-1} \rfloor = \lfloor \frac{\nu n}{n-1} \rfloor = \lfloor \nu + \frac{\nu}{n-1} \rfloor = \nu + \lfloor  \frac{\nu}{n-1} \rfloor =\nu$. Note that $\nu = \frac{m}{n} < n-2$. Therefore,
$\lfloor \frac{m-\nu}{n-2} \rfloor = \lfloor \frac{n\nu-\nu}{n-2} \rfloor = \lfloor \nu + \frac{\nu}{n-2} \rfloor = \nu+ \lfloor  \frac{\nu}{n-2} \rfloor = \nu = \kappa$.
\hfill $\qed$

\vspace{.05in}

We are now in a position to prove Theorem \ref{th:pnas-spectrum}.


\vspace{0.05in}

{\bf Proof of Theorem \ref{th:pnas-spectrum}:} 
We will prove the following claim.

{\bf Claim:} For any $n\ge 2$ and $n-1\le m\le n(n-1)$,
$$p_{L(n,m)}(\lambda) = \lambda(\lambda-\kappa)^{(\kappa+1)(n-1)-m} (\lambda-\kappa-1)^{m-\kappa(n-1)}.$$

Recall that for each pair of $n$ and $m$, distinct $\bbb G(n,m)$ may be generated by Algorithm 1, depending on the choice of the initial tree $\bbb G(n,n-1)$. It is worth emphasizing that we will show the claim holds for all possible $\bbb G(n,m)$, that is, the above characteristic polynomial equation is satisfied for all possible $L(n,m)$.

Note that $1+[(\kappa+1)(n-1)-m] + [m-\kappa(n-1)] = n$. The claim implies that $L(n,m)$ has one eigenvalue at 0, $(\kappa+1)(n-1)-m$ eigenvalues at $\kappa$, and $(\kappa+1)(n-1)-m$ at $\kappa+1$, which together constitute the entire spectrum of $L(n,m)$. The theorem then immediately follows from the claim. Thus, to prove the theorem, it is sufficient to establish the claim. We will prove the claim by induction on $n$.

In the base case when $n=2$, all possible values of $m$ are $1$ and $2$. According to the algorithm description, $\bbb G(2,1)$ contains one arc, $(1,2)$, and $\bbb G(2,2)$ contains two arcs, $(1,2)$ and $(2,1)$. It is straightforward to verify that the claim holds for both $L(2,1)$ and $L(2,2)$. Note that in this special case, $\bbb G(2,1)$ is unique, and so is $\bbb G(2,2)$. 

For the inductive step, suppose that the claim holds for $n=q$ for all possible values of $m$ in $\{q-1,\ldots,q(q-1)\}$ and all possible $\bbb G(q,q-1)$, where $q\ge 2$ is an integer. 
Let $n=q+1$. Take into account all possible $\bbb G(q+1,q)$ and the corresponding values of $m$, which range from $q$ to $(q+1)q$. 

We first consider the case when $m\in\{q,\ldots,q^2-1\}$.
From Lemma \ref{lem:Lentries}, with $n$ replaced by $q+1$, $[L(q+1,m)]_{i(q+1)}=0$ for all $i\in\{1,\ldots,q\}$. That is, all entries in the $(q+1)$th column of $L(q+1,m)$, except for the last entry $[L(q+1,m)]_{(q+1)(q+1)}$, are zero. The same holds for the $(q+1)$th column of $\lambda I-L(q+1,m)$, whose last entry is equal to $\lambda - [L(q+1,m)]_{(q+1)(q+1)}$. 
Then, the Laplace expansion along the $n$th column of $L(q+1,m)$ yields
\begin{align}
p_{L(q+1,m)}(\lambda) &= \det\big(\lambda I - L(q+1,m)\big) \nonumber\\
&= \big(\lambda - [L(q+1,m)]_{(q+1)(q+1)}\big) p_M(\lambda),\label{eq:characteristic}
\end{align}
where $M$ is the $q\times q$ submatrix of $L(q+1,m)$ obtained by removing the $(q+1)$th row and $(q+1)$th column of $L(q+1,m)$. 
Since all row sums of $L(q+1,m)$ are zero, and all entries in its $(q+1)$th column, except for the last entry, are zero, $M$ also has all row sums equal to zero. It follows that $M$ is the Laplacian matrix of a certain graph $\bbb H$ with $q$ vertices, and $\bbb H$ is the subgraph of $\bbb G(q+1,m)$ induced by the vertex subset $\{1,\ldots,q\}$. 
From Lemma \ref{lm:delete-vertex-n}, with $n$ replaced by $q+1$, $\bbb H=\bbb G(q, m-d_{q+1})$, where $d_{q+1}$ denotes the in-degree of vertex $q+1$ in $\bbb G(q+1,m)$, and 
$\bbb G(q, m-d_{q+1})$ is the graph constructed by Algorithm~1 starting from the subgraph of $\bbb G(q+1,q)$ induced by the vertex subset $\{1,\ldots,q\}$. 
Thus, $M$ is the Laplacian matrix of $\bbb G(q, m-d_{q+1})$. 
Let $\gamma \dfb \lfloor \frac{m}{q+1} \rfloor$. We consider the following two cases separately.

{\bf Case 1:} Suppose that $q+1$ divides $m$, which implies $m = \gamma (q+1)$. From the definition of a Laplacian matrix and Proposition \ref{prop:Galmostregular}, $[L(q+1,m)]_{(q+1)(q+1)} = d_{q+1} = \gamma$. Then, $M$ is the Laplacian matrix of $\bbb G(q, m-d_{q+1})=\bbb G(q,m-\gamma)$. From \eqref{eq:characteristic} and the induction hypothesis,  
\eq{
p_{L(q+1,m)}(\lambda) 
= (\lambda - \gamma) p_M(\lambda),\label{eq:case1a}
}
$$p_M(\lambda) = \lambda(\lambda-\beta)^{(\beta+1)(q-1)-m+\gamma} (\lambda-\beta-1)^{m-\gamma-\beta(q-1)},$$
where $\beta \dfb \lfloor\frac{m-\gamma}{q-1}\rfloor$.
The analysis is further divided into two scenarios based on the value of $m$. First, consider when $m = q^2-1$, which implies $\gamma=q-1$ and thus $\beta=q$. Then, $(\beta+1)(q-1)-m+\gamma = q-1$ and $m-\gamma-\beta(q-1) = 0$. 
It follows from \eqref{eq:case1a} that
\eq{
p_{L(q+1,m)}(\lambda) = \lambda(\lambda-q+1)(\lambda -1)^{q-1}. \nonumber 
}
Meanwhile, $(\gamma+1)q-m = 1$ and $m-\gamma q = q-1$. Thus, the above equation validates the claim with $n$ replaced by $q+1$. 
Next, consider when $1\le m \le q^2-2$. From Lemma~\ref{lm:n-divides-m}, with $n$ replaced by $q+1$, $\gamma = \beta =\alpha \dfb \lfloor\frac{m}{q}\rfloor$. Then, from~\eqref{eq:case1a},
\eq{
p_{L(q+1,m)}(\lambda) = \lambda(\lambda-\alpha)^{(\alpha+1)q-m}(\lambda -\alpha-1)^{m-\alpha q},\label{eq:consistentclaim}
}
which proves the claim with $n$ replaced by $q+1$.

{\bf Case 2:} Suppose that $q+1$ does not divide $m$, which implies $m -\gamma (q+1)>0$. 
From the definition of a Laplacian matrix and Proposition \ref{prop:Galmostregular}, $[L(q+1,m)]_{(q+1)(q+1)} = d_{q+1} = \gamma +1$. Then, $M$ is the Laplacian matrix of $\bbb G(q, m-d_{q+1})=\bbb G(q,m-\gamma-1)$. From \eqref{eq:characteristic} and the induction hypothesis,  
\eq{
p_{L(q+1,m)}(\lambda) 
= (\lambda - \gamma-1) p_M(\lambda),\label{eq:case1b}
}
$$p_M(\lambda) \! = \! \lambda(\lambda-\beta')^{(\beta'+1)(q-1)-m+\gamma'} (\lambda-\beta'-1)^{m-\gamma'-\beta'(q-1)},$$
where $\gamma'=\gamma+1$ and $\beta' \dfb \lfloor\frac{m-\gamma-1}{q-1}\rfloor$.
With $n$ replaced by $q+1$, Lemma \ref{lm:floor} and Lemma \ref{lm:p-k} respectively imply that $\beta'=\alpha$ and $\gamma\in\{\alpha-1,\alpha\}$. 
The analysis is then divided into two scenarios based on the value of $\gamma$. 
First, suppose $\gamma = \alpha-1$. Then, $(\beta'+1)(q-1)-m+\gamma' = (\alpha+1)(q-1)-m+\alpha$ and $m-\gamma'-\beta'(q-1) = m-\alpha q$. It follows that \eqref{eq:case1b} simplifies to \eqref{eq:consistentclaim}, 
which validates the claim.  
Next, suppose $\gamma = \alpha$. 
Then, $(\beta'+1)(q-1)-m+\gamma' = (\alpha+1)q-m$ and $m-\gamma'-\beta'(q-1) = m-\alpha q-1$. With these equalities, \eqref{eq:case1b} once again leads to \eqref{eq:consistentclaim}, thereby proving the claim.

The two cases above collectively establish the inductive step for $m\in\{q,\ldots,q^2-1\}$.
In what follows, we address the remaining case where $m \in\{q^2,\ldots, q^2+q\}$. 
From Lemma \ref{prop:G-large-m}, with $n$ replaced by $q+1$, $\overline{\bbb G(q+1, m)}$, the complement of $\bbb G(q+1,m)$, is the $(q+1)$-vertex directed forest consisting of a directed star with $(q^2+q-m+1)$ vertices, rooted at vertex $q+1$, and $m-q^2$ isolated vertices. 
Thus, $\overline{\bbb G(q+1,m)}$ is acyclic, with $(q^2+q-m)$ vertices of in-degree $1$ and the other $(m-q^2+1)$ vertices of in-degree $0$. 
From Lemma \ref{lm:acyclic}, the Laplacian spectrum of $\overline{\bbb G(q+1,m)}$ consists of $(q^2+q-m)$ eigenvalues at $1$ and $(m-q^2+1)$ eigenvalues at $0$. 
Then, from Lemma \ref{lm:complement}, the Laplacian spectrum of $\bbb G(q+1,m)$ is composed of a single eigenvalue of $0$, an eigenvalue of $q$ with multiplicity $q^2+q-m$, and an eigenvalue of $q+1$ with multiplicity $m-q^2$. 
This leads to 
$$p_{L(q+1,m)}(\lambda) = \lambda (\lambda-q)^{q^2+q-m}(\lambda-q-1)^{m-q^2},$$
which validates the claim with $n$ replaced by $q+1$. We therefore complete the proof of the inductive step.
\hfill $\qed$

\section{Discussion}

We have resolved a long-standing conjecture on optimal directed graph topologies for network synchronization proposed in \cite{Nishikawa10}. 
The conjecture states that the normalized Laplacian eigenvalue spread is minimized when the Laplacian spectrum follows a specific pattern. A minimal normalized Laplacian eigenvalue spread indicates optimal synchronizability of the network, as demonstrated in simulations in \cite{Nishikawa10}. 
We establish the conjecture in a stronger form (cf. Theorem \ref{thm:main}) by not only proving that the conjectured Laplacian spectrum pattern is a necessary and sufficient algebraic condition, but also showing that this condition is always realizable. In particular, we prove that the Laplacian spectrum pattern is achievable for any feasible pair of vertex and arc numbers by a class of almost regular directed graphs, which can be systematically and efficiently constructed via an inductive algorithm.

The algorithm proposed in this paper is motivated by our recent work \cite{acc25}, where an algorithm was designed to achieve fast/fastest consensus while its generated graphs also possess the Laplacian spectrum described in \eqref{eq:spectrum}. It turns out the algorithm here subsumes the algorithm in \cite{acc25} as a special case, as shown by the following lemma.  



\vspace{.05in}

\begin{lemma}\label{lm:acc-is-special-case}
If $\mathbb{G}(n, n-1)$ is a directed star, then for all $n \geq 2$ and $n-1 \leq m \leq n(n-1)$, the graph $\mathbb{G}(n, m)$ is identical to the one constructed by the algorithm in \cite{acc25} with $n$ vertices and $m$ arcs.
\end{lemma}

\vspace{.05in}

{\bf Proof of Lemma \ref{lm:acc-is-special-case}:}
If $\mathbb{G}(n, n-1)$ in Algorithm 1 is set to be a directed star, it is rooted at vertex $1$. It is easy to verify that when $m = n-1$, the graph constructed by the algorithm in \cite{acc25} is also a directed star rooted at vertex~$1$ (cf. Proposition 1 in \cite{acc25}). 
Thus, to prove the lemma, it remains to consider the case $m\ge n$. 
From Proposition \ref{prop:Galmostregular}, the graph $\mathbb{G}(n, m)$ constructed by Algorithm 1 has the in-degree sequence given in \eqref{eq:indegree}. We use Lemma \ref{lem:incoming-neighbors} to determine the in-neighbors of each vertex. The $d_1$ incoming arcs of vertex $1$ originate from $d_1$ vertices whose indices are in $\{2,\ldots,d_1+1\}$. For each $i\in\{2,\ldots,n\}$, the $d_i$ incoming arcs of vertex $i$ originate from vertex $i_j$, the unique vertex index such that $(i_j,i)$ is an arc in $\bbb G(n,n-1)$, and from $d_i-1$ other vertices whose indices are the $d_i-1$ smallest elements of $\{1,\ldots,n\} \setminus \{i, i_j\}$. Since $\mathbb{G}(n, n-1)$ is a directed star rooted at vertex $1$, $i_j = 1$ for all $i\in\{2,\ldots,n\}$. Then, the in-neighbor indices for vertex $i$ simplify to the $d_i$ smallest elements of $\{1, \ldots, n\} \setminus \{i\}$. 
It is easy to see that this characterization of the in-neighbor set holds for all vertices including vertex $1$.
It has been shown in \cite[Proposition 2 and Lemma 6]{acc25} that the graph with $n$ vertices and $m$ arcs, constructed by the algorithm presented there, has the same in-degree sequence and the same in-neighbor set characterization for each vertex. 
Therefore, for all $n \geq 2$ and $n-1 \leq m \leq n(n-1)$, the two graphs constructed respectively by Algorithm 1 starting with a directed star and by the algorithm in \cite{acc25} are always identical.
\hfill $\qed$

\vspace{.05in}

\begin{figure}[h!]
\centering
\includegraphics[width=2.2in]{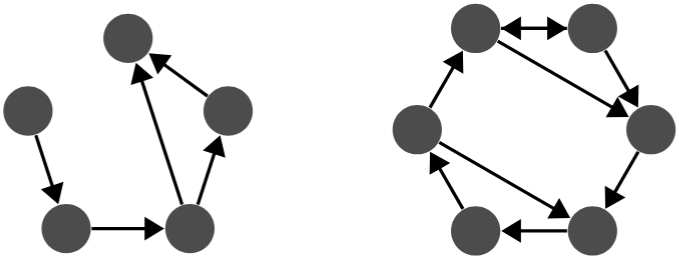} 
\caption{Two optimal graphs with a minimal Laplacian eigenvalue spread that are not generated by Algorithm 1}
\label{fig:notconstructed}
\end{figure}

All graphs generated by Algorithm 1 are almost regular and achieve the conjectured minimum (normalized) Laplacian eigenvalue spread. There exist optimal graphs with minimal eigenvalue spread that are not constructed by Algorithm 1. Moreover, optimal graphs are not necessarily almost regular. Figure \ref{fig:notconstructed} illustrates these.
Specifically, the left graph in Figure \ref{fig:notconstructed} is an optimal graph with 5 vertices and 5 arcs. The left graph is not almost regular, while the graphs constructed by Algorithm 1, $\bbb G(5,5)$, are almost regular. The right graph is an optimal graph with 6 vertices and 9 arcs, in which every vertex has out-degree at least 1. In contrast, in $\bbb G(6,9)$, vertex~6 always has out-degree 0 (cf. Lemma \ref{lem:Lentries}).

It remains unclear how to identify all directed graphs that achieve the minimal Laplacian eigenvalue spread. Even when restricting attention to the class of almost regular directed graphs, a complete characterization of all optimal topologies is still lacking. 
Addressing these open problems and developing systematic methods to enumerate or recognize such optimal graphs are natural directions for future research, albeit challenging ones.

All directed graphs constructed by Algorithm~1 are not only optimal for network synchronization, but also typically guarantee fast convergence for continuous-time linear consensus processes.
A continuous-time linear consensus process over a simple directed graph $\bbb G$ is modeled by a linear differential equation of the form $\dot x(t)=-Lx(t)$, where $x(t)$ is a vector in $\R^n$ and $L$ is the Laplacian matrix of $\bbb G$ \cite{reza1}.
From standard linear systems, the system reaches a consensus as $t \to \infty$, that is, $\lim_{t \to \infty} x(t) = a\mathbf{1}$ for some constant $a$, exponentially fast if and only if the second smallest real part among all eigenvalues of $L$ is positive. This spectral quantity determines the worst-case convergence rate of the consensus process; accordingly, we refer to it as the algebraic connectivity of directed graph $\mathbb{G}$ and denote it by $a(\mathbb{G})$.

\vspace{.05in}

\begin{lemma}\label{lem:aG-directed-upper}
For any simple directed graph $\bbb G$ with $n$ vertices and $m$ arcs, $a(\bbb G) \le \frac{m}{n-1}$.
\end{lemma}

\vspace{.05in}

{\bf Proof of Lemma \ref{lem:aG-directed-upper}:}
Label the $n$ Laplacian eigenvalues of $\bbb G$ as $\lambda_1, \lambda_2,\ldots, \lambda_n$ such that ${\rm Re}(\lambda_1)\le {\rm Re}(\lambda_2)\le \cdots \le {\rm Re}(\lambda_n)$. It is clear that $\lambda_1 = 0$.
Note that the $i$th diagonal entry of the Laplacian matrix equals the in-degree of vertex $i$, denoted by $d_i$. 
Since the sum of all eigenvalues of a matrix equals its trace, 
$\sum_{i=1}^n \lambda_i = {\rm tr}(L(\bbb G)) = \sum_{i=1}^n d_i = m$. Thus,
$ m = \sum_{i=2}^n {\rm Re}(\lambda_i) \ge (n-1){\rm Re}(\lambda_2) = (n-1)a(\bbb G)$, 
which implies $a(\bbb G) \le \frac{m}{n-1}$.
\hfill $\qed$


From Theorem~\ref{th:pnas-spectrum}, all directed graphs constructed by Algorithm 1 with $n$ vertices and $m\ge n-1$ arcs have algebraic connectivity equal to $\lfloor \frac{m}{n-1} \rfloor$. It follows from Lemma~\ref{lem:aG-directed-upper} that the algebraic connectivity of graphs constructed by Algorithm~1 is generally ``close to'' the maximum possible value, with a gap of less than $1$. Consequently, these directed graphs, all of which are almost regular, typically guarantee fast consensus performance for continuous-time consensus processes implemented on them. It is worth emphasizing that this ``close-to-maximum'' property does not hold for all almost regular directed graphs. Further related results and a more detailed discussion can be found in \cite{acc25}.

We conclude the paper by presenting a couple of complementary results and discussing promising future directions.

\subsection{Minimal $\sigma^2$\! Implies Rootedness}\label{subsec:rooted}

The necessary and sufficient spectral condition~\eqref{eq:spectrum} for optimal network synchronization, that is, minimal $\sigma^2$, is an algebraic condition. At the same time, it is well known that network synchronization requires a certain level of connectivity among the agents. The following lemma shows that this required graphical connectivity is implicitly implied by the algebraic condition.

\vspace{.05in}

\begin{lemma}\label{lem:sigma4rooted}
    If a directed graph $\mathbb{G}$ with $n$ vertices and $m \ge n-1$ arcs attains the minimal value of $\sigma^2$, then $\mathbb{G}$ is rooted. 
\end{lemma}

\vspace{.05in}

Corollary~\ref{coro:pnas} guarantees that $a(\mathbb{G})$ is a positive integer when $m \ge n-1$. Hence, Lemma~\ref{lem:sigma4rooted} is an immediate consequence of the following result.

\vspace{.05in}

\begin{lemma}\label{lm:positive}
For any directed graph $\bbb G$, 
    $a(\bbb G)$ is positive if, and only if, $\bbb G$ is rooted. 
\end{lemma}

\vspace{.05in}

{\bf Proof of Lemma \ref{lm:positive}:}
Let $D_{{\rm out}}$ be the out-degree matrix of $\bbb G$, which is assumed to have $n$ vertices; this matrix is an $n\times n$ diagonal matrix whose $i$th diagonal entry equals the out-degree of vertex $i$. The out-degree based Laplacian matrix is $L_{{\rm out}}=D_{{\rm out}}-A'$, where $A$ is the adjacency matrix of $\bbb G$. 
Recall that $L(\bbb G)= L_{{\rm out}}(\bbb G')$. 
Note that the set of all possible simple directed graphs with $n$ vertices is invariant under the graph transpose operation. 
Lemma 2 in \cite{Wu05} shows that the second smallest real part of all eigenvalues of $L_{{\rm out}}(\bbb G)$ is positive if and only if $\bbb G'$ is rooted, which consequently implies that $a(\bbb G)$ is positive if and only if $\bbb G$ is rooted.
\hfill$\qed$


\subsection{Graphs with Integer Weights}

Corollary~\ref{coro:pnas} is stated for unweighted directed graphs and can be generalized to directed graphs with arbitrary integer weights, including negative integer weights.
To this end, we first introduce the definition of the Laplacian matrix for weighted graphs.

A simple weighted directed graph is a simple directed graph in which each arc is assigned a nonzero real-valued weight.
For any simple weighted directed graph $\bbb G_{{\rm w}}$ with $n$ vertices, let $w_{ij}$ denote the weight of arc $(j,i)$, if it exists. We use $D_{{\rm w}}$ and $A_{{\rm w}}$ to denote the corresponding in-degree matrix and adjacency matrix, respectively.
Specifically, $D_{{\rm w}}$ is an $n\times n$ diagonal matrix whose $i$th diagonal entry equals $\sum_{j=1}^n w_{ij}$, the weighted in-degree of vertex $i$, with the understanding that $w_{ij}=0$ if the arc $(j,i)$ does not exist; $A_{{\rm w}}$ is an $n\times n$ matrix whose $ij$th entry equals $w_{ij}$ if $(j,i)$ is an arc  in $\bbb G_{{\rm w}}$, and equals $0$ otherwise. The Laplacian matrix of $\bbb G_{{\rm w}}$ is defined as $L_{{\rm w}}=D_{{\rm w}}-A_{{\rm w}}$. 

It is easy to see that a weighted Laplacian matrix $L_{{\rm w}}$ always has an eigenvalue at 0 since all its row sums equal 0. 
Since the sum of all eigenvalues of a matrix equals its trace, the sum of all eigenvalues of a weighted Laplacian matrix $L_{{\rm w}}$ equals $m_{{\rm w}} \dfb \sum_{i=1}^n\sum_{j=1}^n w_{ij}$. 
This quantity equals the total sum of the weights of all arcs in the corresponding weighted directed graph $\mathbb{G}_{{\rm w}}$ and is called the net number of arcs\footnote{
The quantity is termed the net number of links in \cite{Nishikawa10}.
} 
of $\mathbb{G}_{{\rm w}}$.
Let $\lambda_1,\lambda_2,\ldots,\lambda_n$ be the $n$ eigenvalues of $L_{{\rm w}}$, with $\lambda_1=0$ and $\lambda_2,\ldots,\lambda_n$ possibly complex, and define the normalized eigenvalue spread of $L_{\mathrm w}$ in the same manner as in \eqref{eq:sigma}.

For integer-weighted directed graphs, that is, when all $w_{ij}$ are integers, $m_{\mathrm w}$ is also an integer. In this case, using the same arguments as those in the proof of Corollary~1, as given in the paragraph immediately preceding its statement, we obtain the following result.

\vspace{.05in}

\begin{corollary}\label{coro:integerweights}
For any simple integer-weighted directed graph with $n$ vertices and net number of arcs $m_{{\rm w}}$, 
$$\sigma^2 \ge \textstyle\frac{1}{(n-1)^2} \big[m_{{\rm w}}-(n-1)\kappa_{{\rm w}}\big]\big[(n-1)(\kappa_{{\rm w}}+1)-m_{{\rm w}}\big],$$ 
and equality holds 
if, and only if, the Laplacian spectrum is
\eq{
0, \underbrace{\kappa_{{\rm w}},\; \ldots, \;\kappa_{{\rm w}},}_{(n-1)(\kappa_{{\rm w}}+1)-m_{{\rm w}}} \; \underbrace{\kappa_{{\rm w}}+1, \ldots, \kappa_{{\rm w}}+1}_{m_{{\rm w}}-(n-1)\kappa_{{\rm w}}}, \label{eq:spectrum4weighted}
}
where $\kappa_{{\rm w}} \dfb \lfloor \frac{m_{{\rm w}}}{n-1}\rfloor$.
\end{corollary}

\vspace{.05in}

It is worth emphasizing that Corollary \ref{coro:integerweights} has the same form as Corollary~\ref{coro:pnas}, except that the number of arcs $m$ is replaced by the net number of arcs $m_{\mathrm w}$, even though individual arc weights may be negative integers and, as a consequence, $m_{\mathrm w}$ may be zero or negative.
This is consistent with the simulation results reported in \cite[Page 10345]{Nishikawa10}, which suggest that the conjecture also applies to simple integer-weighted directed graphs.
Notwithstanding this, Corollary \ref{coro:integerweights} itself does not guarantee that the minimal possible value of $\sigma^2$ is attainable.

In the case when $n-1\le m_{{\rm w}}\le n(n-1)$, the minimal value of $\sigma^2$ can be achieved by considering only unweighted directed graphs, for which $m_{{\rm w}}$ equals the number of arcs, and by constructing corresponding optimal graph(s) via Algorithm~1. 
However, it remains theoretically unclear whether, and under what conditions, the minimal value of $\sigma^2$ can be attained by integer-weighted directed graphs, especially when negative integer weights are involved, although the existence of such graphs has often been observed in simulations in~\cite{Nishikawa10}.

A more interesting case arises when $m_{{\rm w}} < n-1$ or $m_{{\rm w}} > n(n-1)$, in which the optimal directed graphs cannot be unweighted. In particular, when $m_{{\rm w}} < n-1$, negative weights must be involved, since if all weights are positive integers, the graph cannot be rooted, and thus network synchronization cannot be achieved.

Analyzing network synchronization involving negative-weighted interactions poses an important and challenging research direction.
Simulations in \cite{Nishikawa10} have demonstrated that negative-weighted interactions can enhance network synchronization performance, while the Laplacian matrices of weighted graphs with negative weights do not possess the favorable properties of Laplacian matrices of positive-weighted graphs.
In particular, the spectral condition \eqref{eq:spectrum4weighted} does not guarantee rooted connectivity, nor does it ensure that all eigenvalues other than $\lambda_1$ have positive real parts, both of which are necessary conditions for network synchronization. In contrast, for positive-weighted directed graphs, both rooted connectivity and positivity of the real parts of all nonzero eigenvalues are guaranteed under the same spectral condition.
Moreover, for weighted graphs with negative weights, even if the underlying graph is connected, the eigenvalues other than $\lambda_1$ do not necessarily have positive real parts \cite{tac_laplacian}.  
This gap between empirical observations and available spectral theory highlights the need for new analytical tools for network synchronization with negative-weighted interactions.

The above three paragraphs outline a few directions for future theoretical development. In this context, Corollary \ref{coro:integerweights} serves as a useful starting point and a basis for further analysis.

\section*{Acknowledgment}

The work of J. Urschel is supported by the National Science Foundation under grant no. DMS-2513687. 
J. Liu wishes to thank Wei Chen (Peking University) for introducing the conjecture in \cite{Nishikawa10} and Dan Wang (Nanjing University) for discussing an alternative proof in the special tree case, both several years ago.

\bibliographystyle{unsrt}
\bibliography{susie,jicareer}

\end{document}